\documentclass[a4paper,twoside,11pt]{article}
\usepackage{latexsym,amssymb,enumerate,amsmath,epsfig,amsthm,graphicx}
\usepackage[margin=1in]{geometry}
\usepackage{setspace,color}
\usepackage{tikz,url}
\usepackage{floatrow}
\usepackage{subfig}
\usepackage[ruled]{algorithm2e}
\usepackage{csquotes}
\usepackage{dsfont}
\usepackage[multidot]{grffile}
\usepackage{enumerate}
\usepackage{hyperref} 
\hypersetup{
    colorlinks=true,
    linkcolor=blue,
    filecolor=magenta,      
    urlcolor=cyan,
    citecolor=red,
    pdfpagemode=FullScreen,
    }

\newcommand{\x}{\mathbf{x}}

\newcommand{\reminder}[1]{{#1}}

\begin{document}

\title{Within-Cluster Variability Exponent for Identifying Coherent Structures in Dynamical Systems}
\author{
Wai Ming Chau\thanks{Department of Mathematics, the Hong Kong University of Science and Technology, Clear Water Bay, Hong Kong. Email: {\bf wmchau@connect.ust.hk}}
\and
Shingyu Leung\thanks{Department of Mathematics, the Hong Kong University of Science and Technology, Clear Water Bay, Hong Kong. Email: {\bf masyleung@ust.hk}}
}

\date{}

\maketitle

\begin{abstract}
We propose a clustering-based approach for identifying coherent flow structures in continuous dynamical systems. We first treat a particle trajectory over a finite time interval as a high-dimensional data point and then cluster these data from different initial locations into groups. The method then uses the normalized standard deviation or mean absolute deviation to quantify the deformation. Unlike the usual finite-time Lyapunov exponent (FTLE), the proposed algorithm considers the complete traveling history of the particles. We also suggest two extensions of the method. To improve the computational efficiency, we develop an adaptive approach that constructs different subsamples of the whole particle trajectory based on a finite time interval. To start the computation in parallel to the flow trajectory data collection, we also develop an on-the-fly approach to improve the solution as we continue to provide more measurements for the algorithm. The method can efficiently compute the WCVE over a different time interval by modifying the available data points.
\end{abstract}

\section{Introduction}
\label{sec:intro}

Finite-time Lyapunov exponent (FTLE)  \cite{halyua00,hal01,hal01b,shalekmar05,lekshamar07} is a widely used Lagrangian quantity to hint the location of any Lagrangian coherent structure (LCS) in a given velocity field $\mathbf{u}(\mathbf{x}(t),t)$. It measures the rate of change in the distance between neighboring particles across a finite interval of time with an infinitesimal perturbation in the initial position. To obtain the FTLE field, one needs to first compute the flow map which links the initial location of a particle with the arrival position based on the characteristic line. Mathematically these particles in the extended phase space satisfy the ordinary differential equation (ODE) given by
\begin{equation}
\dot{\mathbf{x}}(t) = \mathbf{u}(\mathbf{x}(t),t)
\label{Eqn:RayTracing}
\end{equation}
with the initial condition $\mathbf{x}(t_0)=\mathbf{x}_0$ and a Lipschitz velocity field $\mathbf{u}:\mathbb{R}^d \times \mathbb{R} \rightarrow \mathbb{R}^d$. The flow map $\Phi_{t_0}^{t_0+T}:\mathbb{R}^d \rightarrow \mathbb{R}^d$ is defined as the mapping which takes the point $\mathbf{x}_0$ to the particle location at the final time $t=t_0+T$, i.e $\Phi_{t_0}^{t_0+T}(\mathbf{x}_0)=\mathbf{x}(t_0+T)$ with $\mathbf{x}(t)$ satisfying (\ref{Eqn:RayTracing}). The FTLE is then defined using the largest eigenvalue of the deformation matrix based on the Jacobian of this resulting flow map. In a series of recent studies \cite{leu11,leu13,youwonleu17,youleu18,youleu18b,lywn19,ngleu19,youleu20}, we have developed various Eulerian approaches to numerically compute the FTLE on a fixed Cartesian mesh. The idea is to combine the approach with the level set method \cite{oshset88,oshfed03} which allows the flow map to satisfy a Liouville equation. Such hyperbolic partial differential equations (PDEs) can then be solved by any well-established robust and high order accurate numerical methods. We have also developed several other quantities to extract coherent structures in a given velocity field. But we will not compare them with the proposed approach in this work but concentrate only on the FTLE as a demonstration. We refer interested readers to \cite{youleu21,youleu21b} and references therein for more discussions.

\begin{figure}[!htb]
(a)\includegraphics[width=0.45\textwidth]{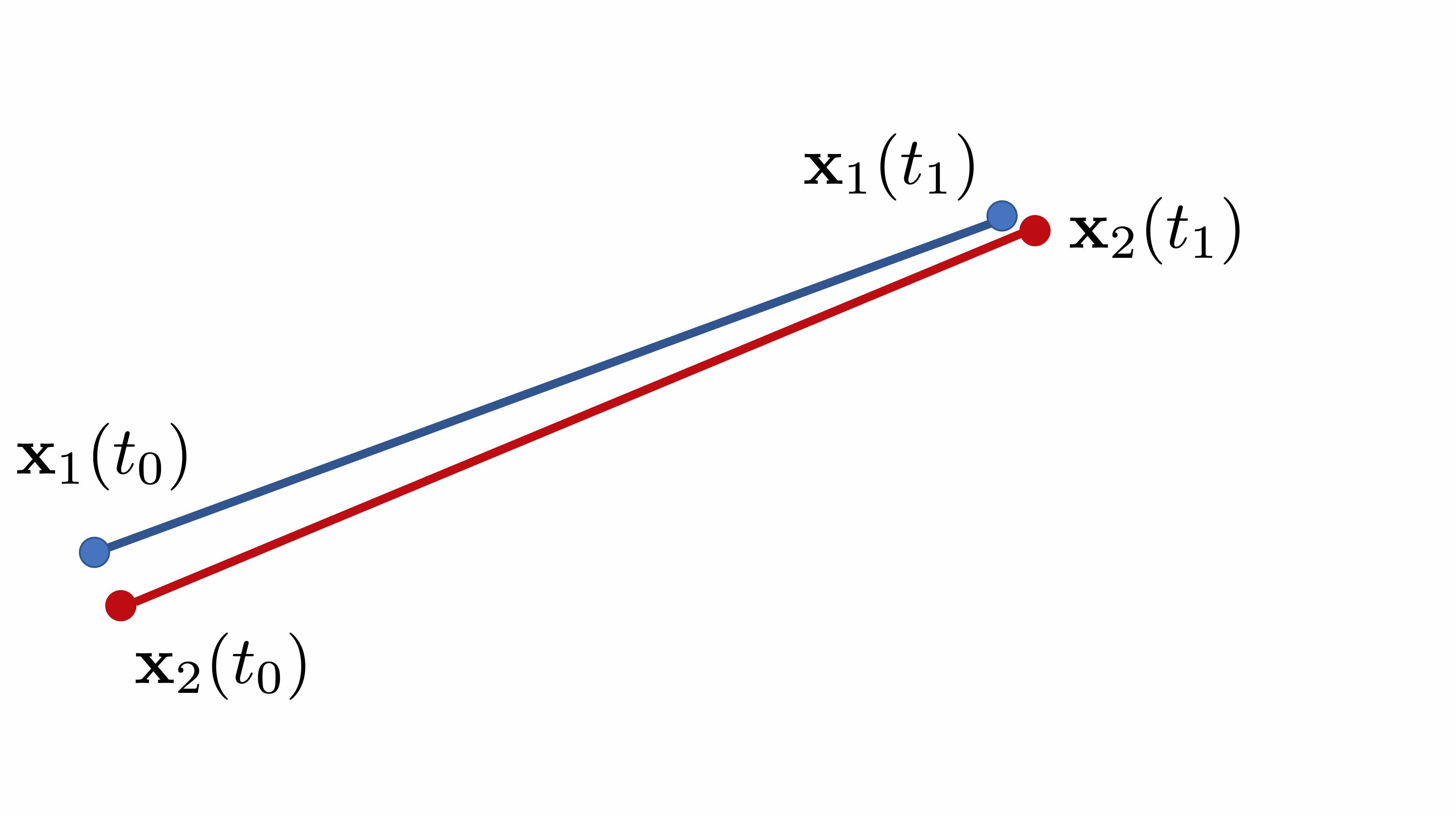}
(b)\includegraphics[width=0.45\textwidth]{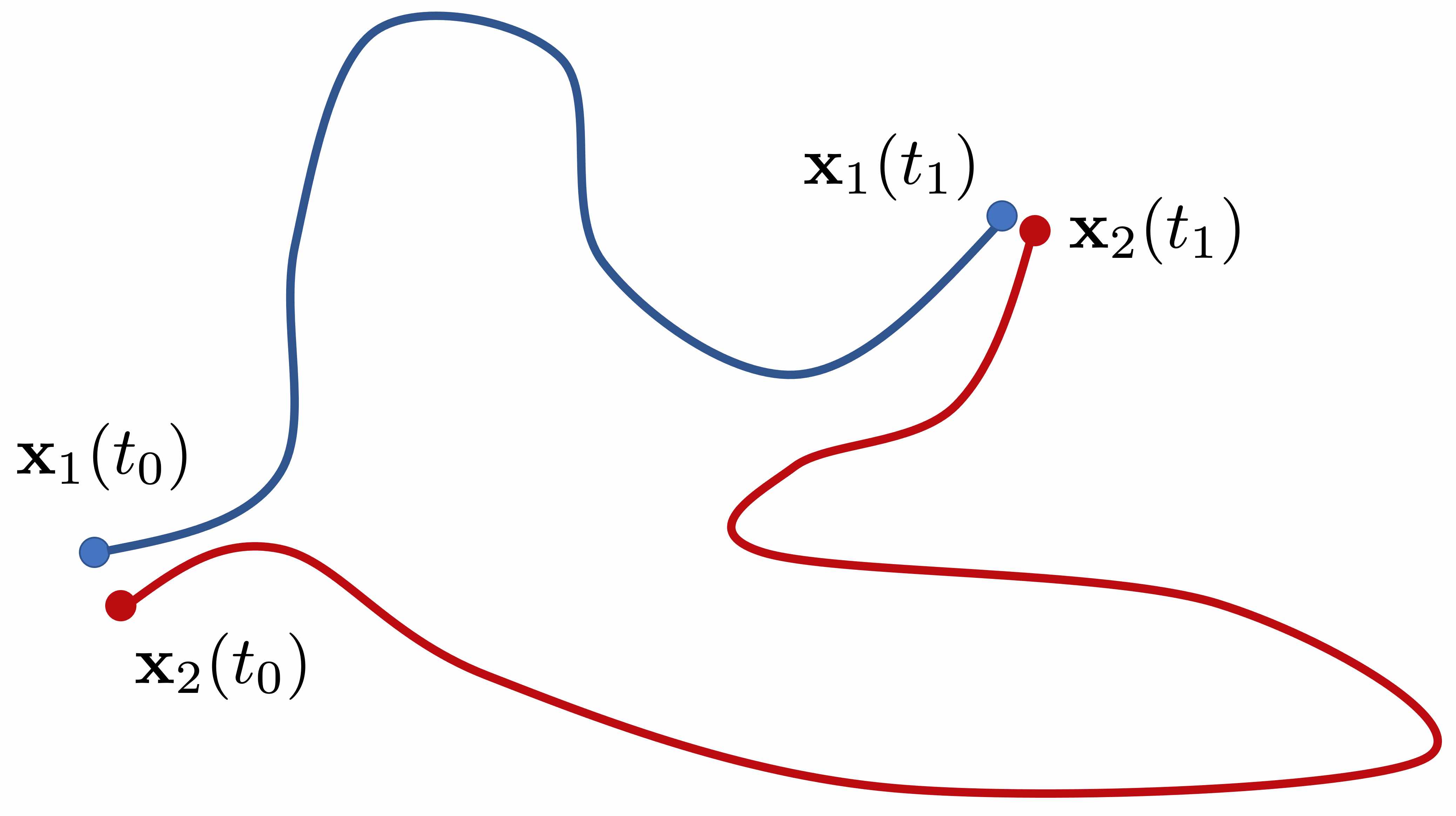}
\caption{FTLE requires only the flow map from the initial location at the initial time to the final arrival location at the terminal time but ignores the particle's intermediate locations.}
\label{Fig:FTLE}
\end{figure}

From the computation of the FTLE, we notice that FTLE ignores the particle's intermediate locations, which might lead to a misinterpretation of coherent structure in a dynamic system. Because the FTLE requires only the flow map from the initial location at the initial time to the final arrival location at the terminal time, the quantity cares only about the particle's initial and final positions. If an initial patch of particles stays close for a significant period but diverts quickly right before the terminal time, FTLE would still give a large quantity. Suppose particles in an initial patch travel very differently with unrelated traveling history but reach a similar region at the terminal time. As shown in Figure \ref{Fig:FTLE}, FTLE ignores any significant dispersion at the intermediate time and returns a small quantity. When neighboring particles arrive at the exact locations, all intermediate locations are irrelevant to the FTLE completely. Therefore, the FTLE for both cases in Figure \ref{Fig:FTLE}(a) and (b) could be the same even though these particles could have two very different paths.

The second issue about the approach is that the FTLE could depend heavily on the time period $T$ chosen for visualization. Taking the example in Figure \ref{Fig:FTLE}(b) again, we will obtain a significantly larger FTLE if we take approximately half of the time interval when these particles are widely separated. This phenomenon could be far from satisfactory in real applications. Even if we obtain a significant FTLE region indicating a large rate of infinitesimal separation, the structure can disappear immediately if we allow the dynamical system to evolve a little longer.

In this paper, we develop an approach to incorporate the complete particle's trajectory and introduce a local measure to quantify the variations in the traveling history within a small neighborhood of each particle. Trajectory clustering is not a new topic. For example, a partition-and-group framework has been developed in \cite{leehanwha07} to discover common sub-trajectories. \cite{vlagunkol02} has developed a non-metric similarity function based on the longest common subsequence to determine objects that move in a similar pattern. Some model-based clustering techniques \cite{gafsmy99,grscg07,wycm11,fkss13} and a more recently developed deep learning approach \cite{yzzhb17} have also been proposed to classify trajectory data into groups. In \cite{youleu14}, we have also developed a coherent ergodic partition method to separate trajectories into clusters. Instead of clustering to the whole trajectory as a high-dimensional data point, we have proposed integrating and computing the long-time averages of a set of functions along the trajectories. This idea efficiently projects the high-dimensional data to a low-dimensional manifold. Then based on these function averages along each trajectory, we partition these data into groups.

However, all these data analysis methods aim to visualize or construct dictionaries and summaries to reveal hidden trajectory patterns or predict possible future routes. The primary purpose has nothing to do with discovering coherent structures in the underlying flow and dynamical system. In this work, we further introduce a measure to quantify the trajectories discrepancy within each cluster and use it to identify any coherent structures in the flow. We call the quantity within-cluster variability exponent (WCVE), approximating the average separation rate between the adjacent particles within a finite time interval. \reminder{As a by-product of the algorithm, the WCVE computation can also identify nonlocal regions that come together later. From those initial locations that belong to the same clusters, particles will move together in a relatively short time (compared to the period $T$ in the overall computations). Unlike the FTLE field, which provides only the local infinitesimal flow information, our proposed new quantity also hints at the global relationship among different initial particles.}

The rest of the paper is organized as follows. In Section \ref{Sec:Background}, we will briefly introduce the finite-time Lyapunov exponent (FTLE). Section \ref{Sec:ProposedApproach} reviews basic definitions of the clustering problem and defines the WCVE for flow coherent structure extraction. Sections \ref{SubSec:Refinement} and \ref{SubSec:OnTheFly} propose two improvements and extensions of the approach. We will develop an adaptive approach to speed up the $k$-mean clustering algorithm when all complete trajectory data are available in the first place. On the other hand, if we continue to collect trajectory measurements, we also develop an \textit{on-the-fly} approach to update our classification results in time. Finally, when we compute the WCVE for another time interval, we propose a simple strategy to update our computational results efficiently. Some two-dimensional examples will be given in Section \ref{Sec:Examples} to demonstrate the effectiveness of our approach.

%%%%%%%%%%
\section{Background: Finite-time Lyapunov Exponent (FTLE)}
\label{Sec:Background}

Collecting the solutions to ODEs (\ref{Eqn:RayTracing}) for all initial conditions $\x(t_0)=\x_0\in\Omega$ at all time $t\in\mathbb{R}$, we introduce the \textit{flow map} $\Phi_a^b:\Omega \rightarrow \mathbb{R}^d$
such that $\Phi_a^b(\x_0)=\x(b)$ represents the arrival location $\x(b)$ at $t=b$ of the particle trajectory satisfying ODEs (\ref{Eqn:RayTracing}) with the initial condition $\x(a)=\x_0$ at the initial time $t=a$. This implies that the mapping will take a point from $\x(a)$ at $t=a$ to another point $\x(b)$ at $t=b$.

The quantity, which measures the separation rate between adjacent particles over a finite-time interval with an infinitesimal perturbation in the initial location, is called the finite-time Lyapunov exponent (FTLE). For example, consider the initial time $t=t_0=0$ and the final time $t=T$, the change in the initial infinitesimal perturbation is given by
$$
		\delta \x(T) =\Phi^T_0(\x+\delta\x(0))-\Phi^T_0(\x) 
		= \nabla \Phi^T_0(\x)\delta\x(0)+O(\|\delta\x(0)\|^2)\,.
$$
The magnitude of the leading order term of this perturbation is given by
$$
	\|\delta\x(T)\|=\sqrt{\langle \delta\x(0),[\nabla\Phi^T_0(\x)]^{\prime}\nabla\Phi^T_0(\x)\delta\x(0) \rangle}
$$
where $(\cdot)^{\prime}$ denotes the transpose of a matrix. Denote $\varDelta^T_0(\x)=[\nabla\Phi^T_0(\x)]^\ast\nabla\Phi^T_0(\x)$ the Cauchy-Green deformation tensor.
The maximum of $\|\delta\x(T)\|$ is achieved when $\delta\x(0)$ is chosen to be aligned with the eigenvector $\textbf{e}(0)$, which is associated with the maximum eigenvalue of $\varDelta^T_0(\x)$, i.e. $\max_{\delta\x(0)}\|\delta\x(T)\|=\sqrt{\lambda_{\max}(\varDelta^T_0(\x)}\|\textbf{e}(0)\|=e^{\sigma^T_0(\x)|T|}\|\textbf{e}(0)\|$. Hence, the FTLE $\sigma^T_0(\x)$ is defined as
\begin{equation}
\sigma^T_0(\x)=\frac{1}{|T|}\ln\sqrt{\lambda_{\max}(\varDelta^T_0(\x))} \, .
\label{eq:FTLE}
\end{equation}

Usual numerical methods in computing the flow map $\Phi_0^T(\mathbf{x})$ rely on ray tracing by solving the ODE (\ref{Eqn:RayTracing}) or the corresponding level set equation \cite{leu11} on a Cartesian mesh at the initial time $t=t_0=0$. Then for each grid point in the computational domain, one constructs the Cauchy-Green deformation tensor by finite differencing the \textit{forward} flow map $\Phi_0^T(\mathbf{x}^i_g)$ on the Cartesian mesh at $t=0$. The \textit{forward} FTLE can then be computed using equation (\ref{eq:FTLE}).

However, whether we use the Lagrangian (ODE) or the Eulerian (PDE) approach, FTLE requires only the flow map from the initial location at the initial time to the final arrival location at the final time. The quantity does not pay attention to the intermediate trajectory but only to the particle's initial and final positions. If an initial patch of particles stays close after a finite period, FTLE ignores any significant dispersion at the intermediate time and returns a small quantity. This observation is, unfortunately, unsatisfactory. Even though the trajectory of these particles could be significantly distinct, the FTLE completely ignores the information and cannot identify such a unique structure in the flow.

\section{Our Proposed Approach}
\label{Sec:ProposedApproach}

This section first summarizes our trajectory clustering problem and introduces some notations that will be useful for later discussions. Then we will present our proposed within-cluster variability exponent (WCVE). We will also present two extensions of the method based on the availability of the trajectory data. If we already have all complete trajectory data on hand, we develop an adaptive refinement approach to improve the computational speed of solving the $k$-means clustering problem. On the other hand, if we continually collect trajectory measurements, we present an \textit{on-the-fly} approach to update the $k$-means solution. Finally, the last section will comment on the numerical method for obtaining the trajectory data. We will briefly discuss the Eulerian approach and how we incorporate the Eulerian flow map into the computations.

\subsection{The $k$-mean Trajectory Clustering Problem}
\label{Sec:k-means}

For a time interval $[t_0,t_0+T]$, we can define a time step $\triangle t=T/n$, for positive $n$ and time $t_i=t_0+i\triangle t$, where $1\leq i\leq n$. Solving the flow equation (\ref{Eqn:RayTracing}), we can obtain the flow map $\Phi_{t_{i-1}}^{t_i}(\mathbf{x}({t_{i-1}}))=\mathbf{x}(t_i)$ for each $i$. Therefore, we define the trajectory $\mathbf{X}$ as the $d\times (n+1)$ matrix with the initial condition $\mathbf{x}(t_0)=\mathbf{x}_0\in\mathbb{R}^d$. 
\begin{equation}
    \mathbf{X}=\begin{pmatrix}
\mathbf{x}(t_0)&
\mathbf{x}(t_1)&
\ldots&\mathbf{x}(t_n)
\end{pmatrix} \in \mathbb{R}^{d\times (n+1)} \, .
\label{Eqn:Tra}
\end{equation}
Suppose we have $M$ trajectories with different initial conditions, we denote $\mathbf{X}_j$ as the $j$-th trajectory with the initial condition $\mathbf{x}_j(t_0)=x_j$. Clustering aims to partition a set of trajectories $\{\mathbf{X}_j:j=1,2,...,M\}$ into $k$ disjoint clusters $C=\{C_l:l=1,2,...,k\}$. We define $\lambda_j \in \{1, 2, . . . , k\}$ cluster label of $\mathbf{X}_j$ (i.e. $\mathbf{X}_j \in C_{\lambda_j}$) and cluster label vector $\Lambda=\{\lambda_1,\lambda_2,...,\lambda_M\}$ to show the clustering result. 

We first solve the dynamical system (\ref{Eqn:RayTracing}) numerically using the Lagrangian or the Eulerian approach. Then, we partition a set of trajectories $\mathbf{X}_j$ with the initial condition $\mathbf{x}_j(t_0)=x_j$ in the domain $\Omega$ by minimizing the within-cluster sum of squares (WCSS) to cluster $M$ trajectories,
\begin{equation}
\mbox{WCSS}(C)=\sum_{l=1}^{k}\sum_{\mathbf{X}\in C_l}||\mathbf{X}-\mathbf{m}_l||_{2}^{2} 
\label{Eqn:WCSS}
\end{equation}
where $\mathbf{m}_l=\frac{1}{\mbox{card}(C_l)}\sum_{\mathbf{X}\in C_l}\mathbf{X}$ represents the mean of the cluster, the notation $||.||_2$ is the $L_2$ norm operator and $\mbox{card}(\cdot)$ represents the cardinality of a set. This problem can be efficiently solved by the $k$-means algorithm. In general, the $k$-means algorithm iterates between an assignment step and an update step. In the assignment step, the algorithm assigns each data point $\mathbf{X}_j$ to one of the classes $C_l$ by looking for the class that minimizes the distance between $\mathbf{X}_j$ and the cluster centers $\mathbf{m}_l$. Once the assignment step has updated all labellings, the update step will recalculate the centers $\mathbf{m}_l$ within each cluster. The iteration stops when there is no change in the assignment or the update step.

\reminder{Note that the discussion in this section is nothing but just a simple trajectory clustering algorithm, see for example \cite{fropad15,hkth16,bankol17}. We are going to discuss our main contributions in the next section.}

\subsection{Within-Cluster Variability Exponent (WCVE)}

Unlike typical trajectory clustering methods, which mainly try to reveal hidden trajectory patterns or predict possible future routes, we observe that the clustering results might provide intrinsic information on the dynamical system. We are not interested in the trajectory pattern but in locations where adjacent trajectories vary their motion patterns. Near these coherent structures of the dynamical system, a tiny perturbation in the initial location will lead to a significant change in the arrival position at a final time and, more importantly, the trajectory pattern.

Unless all particles travel the same way, we always see some trajectory variations within each cluster when the number of classes is smaller than the number of data points, i.e. $k<M$. When the value of $k$ becomes more significant, it will be easier to isolate those particles near the coherent structure of the flow. Data are closer to each other in the high dimensions for those clusters containing particles away from these coherent structures. However, for those clusters containing particles initially put near these flow structures, we see trajectory data will spread more significantly in the high dimensional space. Therefore, we propose determining each cluster's statistical dispersion and variability as a measure of coherent structure in the dynamical system.

Mathematically, we define the within-cluster variability exponent (WCVE) for each grid location $\mathbf{x}_j$ (i.e. each trajectory initial takeoff location) as
\begin{eqnarray}
\mbox{WCVE}(\mathbf{x}_j;t_0,t_n) &=& 
\log \left[ \frac{1}{(n+1)}
\sqrt{\frac{1}{\mbox{card}\left(C_{\lambda_j} \right)-1} \sum_{\mathbf{X}\in C_{\lambda_j}}||\mathbf{X}-\mathbf{m}_{\lambda_j}||^2_2} \right] \, .
\label{Eqn:WCVE} 
\end{eqnarray}
This quantity first measures the standard deviation (SD) of these trajectory data within the cluster in the $d\times(n+1)$-dimensional space. Then, it normalizes the SD by the length of the data, followed by taking a logarithm. The normalization $n$ comes from the size of the time period. Assume that we are looking at a fixed time interval from $t_0$ to $T$ with the discretization $\Delta t=(T-t_0)/n$. Reducing the timestep size $\Delta t$ by half would imply doubling the length of the vector $\mathbf{X}$. This extra factor of $(n+1)$ provides the correct scaling of the expression. \reminder{Finally, the logarithm operator can capture the exponential deviation of the SD among different clusters. This step is, therefore, an analogy to the FTLE and can highlight regions with a significantly large SD.}

Indeed, it is also possible to replace the SD with the mean absolute deviation (MAD) to obtain 
\begin{eqnarray}
\log \left[ \frac{1}{(n+1)}
\frac{1}{\mbox{card}\left(C_{\lambda_j} \right)} \sum_{\mathbf{X}\in C_{\lambda_j}}||\mathbf{X}-\mathbf{m}_{\lambda_j}||_2 \right] \, .
\label{Eqn:MADE} 
\end{eqnarray}
This quantity first measures the MAD of each trajectory data within the cluster in the $d\times(n+1)$-dimensional space. Then, it normalizes the MAD by the length of the data, followed by taking a logarithm.

Other than the obvious argument based on the geometrical interpretation, there are other apparent reasons why MAD is preferred in some applications. Since particles could spread over the whole domain and small perturbation in the initial location could lead to a significant change in the trajectory pattern, we do not expect these trajectory data in the high dimensional space to satisfy any standard or even normal distribution. The SD does not necessarily have better data efficiency and does not provide a more intrinsic understanding of the data distribution. Moreover, SD might over-emphasize the outlier. In each cluster formed by the $k$-means clustering, the approach based on SD might over-react to those minor variations in the high dimensional data points. Having said that, however, we do not observe much difference between the WVCE solutions based on either SD or MAD.

Therefore, in computations, we prefer SD to MAD in the definition of WCVE. The main reason is that such a definition is more aligned with the WCSS in the $k$-means clustering step. The minimizer from the $k$-means clustering algorithm gives the variance within each cluster. We only need to further process the output from the $k$-means to obtain the required quantity.

Here are some other properties of the WCVE.
\begin{itemize}
\item The choice of the vector norm is not unique. In the current definition, we are using the typical 2-norm. It seems more natural to the application since we have data points in the usual Euclidean space $\mathbb{R}^d$. Indeed, it involves a square root which might introduce a complication in the implementation. In this case, we can replace the definition with the 1-norm so that we can compute the quantity in a component-by-component fashion. 
\item Since the dependence in $\mathbf{x}_j$ is reflected only in $\lambda_j$, WCVE is a constant for all locations $\mathbf{x}_j$ within the same cluster labeled by $\lambda_j$. This property implies that WCVE assigns value to individual clusters instead of trajectory data. In visualization, however,  we assign the WCVE to grid points representing the initial take-off location of all trajectories. It will be more helpful for visualizing the coherent structure of the dynamical system.
\item For a fixed location $\mathbf{x}_j$, the WCVE could be discontinuous in $t_n$ in general. In other words,  incorporating more trajectory measurements can sharply change the WCVE field at a specific location. We have this property since increasing the data dimensions might change the $k$-means clustering results, assigning the trajectory starting at $\mathbf{x}_j$ to a different cluster. 
\item As we increase the value of $k$ (the number of clusters), the WCSS (\ref{Eqn:WCSS}) based on the SD reduces. Now, consider the WCVE defined using the MAD. Since the variance bounds the MAD from above, we expect to see a smaller WCVE in general as we reduce the WCSS. In particular, as $k$ approaches the total number of trajectories (i.e. $M$), the WCSS converges to zero since the centroids of $M$ clusters will precisely equal to $\{\mathbf{X}_1,\mathbf{X}_2,...,\mathbf{X}_M\}$. This result implies that the WCVE converges to zero as $k$ approaches $M$. 
\item If we ignore the computation complexity of the $k$-means algorithm, the complexity of determining the WCVE is $O(dnM)$, independent of the number of clusters $k$. We only need to go over all $M$-trajectory data points once. Since we have already obtained the clustering labels, we can continuously update the within-cluster SD or MAD as we process the trajectories. When we have these $k$ numbers after reading through all data points, we assign these numbers to $M$ initial takeoff locations of these trajectories. 
\end{itemize}

\subsection{An Adaptive Refinement Approach}
\label{SubSec:Refinement}

One bottleneck of the method is the computational time of the $k$-means algorithm when we have a large number of particle trajectories in extremely high dimensions. In this section, we develop an adaptive approach to improve the convergence of the $k$-means algorithm. There are existing algorithms to improve the clustering method for general cases. But in this work, we improve the computation efficiency based on the continuity of each particle trajectory. We decompose the trajectory data (\ref{Eqn:Tra}) by subsampling the solution at different time levels, which effectively lowers the dimension of each data point. The $k$-means clustering method converges quickly to the mean vectors in this lower-dimensional space. Then we gradually fill in the intermediate positions of these trajectories and bring up the dimension of the data points. The clustering result from the low dimensions can provide a good initial guess of the $k$-means algorithm in the high dimensions. We simply take the mean of the high dimensional trajectory data to fill in the missing information in the cluster mean due to the dimension raising. 
The adaptive approach consists of the following four main steps.

\begin{itemize}
\item[] \textbf{Step 1:} We decompose the trajectory data (\ref{Eqn:Tra}) as follows: 
\begin{align*} 
&\mathbf{X}_j^{1}=\begin{pmatrix}
                \mathbf{x}_j(t_0)&
                \mathbf{x}_j(t_2)&
                \mathbf{x}_j(t_4)&
                \ldots&\mathbf{x}_j(t_n)
            \end{pmatrix} \in \mathbb{R}^{d\times \frac{n+1}{2}} \\
&\mathbf{X}_j^{2}=\begin{pmatrix}
                \mathbf{x}_j(t_0)&
                \mathbf{x}_j(t_4)&
                \mathbf{x}_j(t_8)&
                \ldots&\mathbf{x}_j(t_n)
            \end{pmatrix} \in \mathbb{R}^{d\times \frac{n+1}{2^2}} \mathsf{...}\\
&\mathbf{X}_j^{N}=\begin{pmatrix}
                \mathbf{x}_j(t_0)&
                \mathbf{x}_j(t_{2^N})&
                \mathbf{x}_j(t_{2^{N+1}})&
                \ldots&\mathbf{x}_j(t_n)
            \end{pmatrix} \in \mathbb{R}^{d\times \frac{n+1}{2^N}}
\end{align*} 
At each level, we subsample each trajectory data and lower the dimension of each point by half. Therefore, we have in total of $M$ data points, with each having dimensions of $d\times \frac{n+1}{2^N}$ at the final ($N$-th) refinement level.   
\item[] \textbf{Step 2:} Apply the $k$-means algorithm to cluster $\left\{ \mathbf{X}_j^{N} \right\}$ and obtain $k$ clusters $C_l$.
\item[] \textbf{Step 3:} Up-sample these $M$ data. Apply the $k$-means algorithm to cluster $\left\{ \mathbf{X}_j^{N-1} \right\}$ with the initial condition of the mean vectors 
$$
\frac{1}{\mbox{card}(C_l)} \sum_{\lambda_j=l} \mathbf{X}_j^{N-1}
$$
determined using the clustering results $C_l$ obtained in the previous stage.
\end{itemize}
We then iterate the last step in the algorithm until we consider $\left\{ \mathbf{X}_j^1 \right\}$ the data points in the most refined level with each data having the size $d\times (n+1)$.

In \textbf{Step 3}, we take the average of the whole trajectory $\mathbf{X}_j^{N-1}$. But for those sampling locations that we have already considered in the $N$-th level, their average leads to the same mean vector in the previous stage. When the computation efficiency at this initialization step is crucial, we can use the clustering results at the previous stage to compute the missing components in the mean vectors and keep the information we have already obtained.

In practice, we do not need to obtain the precise converged result in the $k$-means clustering at any intermediate level since the overall adaptive refinement algorithm uses those solutions only as an initial condition at the next refined level. We find that this strategy significantly improves the convergence speed of the $k$-means clustering step of our computations. 

\reminder{On the other hand, since the solution to the $k$-means algorithm depends heavily on the data, this proposed adaptive algorithm might lead to a slightly different classification result than working on the entire dataset at once. We will demonstrate this in the example section when we look at the minimized value of WCSS from different algorithms. Even though there is no theoretical analysis to quantify the difference, our numerical simulations have shown that we can capture some very similar structures in the solution with various levels of adaptivity.}

\subsection{An \textit{On-the-fly} Approach}
\label{SubSec:OnTheFly}

We have assumed that the complete trajectory data are available in the first place by solving the ODE system for all initial conditions or by measuring the motion trajectory of all objects for all time. In practice, however, we might want to develop an \textit{on-the-fly} algorithm to update the clustering results in parallel with data collection. This approach then not only improves the time \reminder{and memory} management of the computations but also provides a better initial condition for the $k$-means algorithm. 

We could start the computations immediately when we have sufficient information of the trajectories. We first define the intermediate trajectory $\mathbf{X}_{j}^{z}$, for $0<z<n$, as follows 
 \begin{align*} 
    \mathbf{X}_{j}^{z}=\begin{pmatrix}
                    \mathbf{x}_j(t_0)&
                    \mathbf{x}_j(t_1)&
                    \mathbf{x}_j(t_2)&
                    \ldots&\mathbf{x}_j(t_z)
                    \end{pmatrix} \in \mathbb{R}^{d\times (z+1)}
\end{align*}
Our proposed \textit{on-the-fly} approach consists of the following steps.
\begin{itemize}
\item[] \textbf{Step 1:} Apply the $k$-means algorithm and obtain $k$ mean vectors 
$$
\left\{ \mathbf{m}^{z}_l \in \mathbb{R}^{d\times (z+1)} :l=1,2,...,k \right\}
$$ 
and the cluster label vector $\Lambda=\{\lambda_1,\lambda_2,...,\lambda_M\}$ through clustering $\left\{\mathbf{X}_{j}^{z}\right\}$.
\item[] \textbf{Step 2:} Obtain the latest arrival locations $\mathbf{x}_j(t_{z+1})$ at $t=t_{z+1}$ for each of the trajectories.
\item[] \textbf{Step 3:} Apply the $k$-means algorithm to cluster $\left\{ \mathbf{X}_j^{z+1} \right\}$ with the initial condition of the mean vectors 
$$
\frac{1}{\mbox{card}(C_l)} \sum_{\lambda_j=l} \mathbf{X}_j^{z+1}
$$
determined using the clustering results $C_l$ obtained in the previous stage.
\end{itemize}

\reminder{We denote the above algorithm \textit{on-the-fly-1} with the index equal to 1, representing that the $k$-means algorithm is run for every time step when we have the most updated location of all the data points. In terms of computational complexity, even though the classification results from the previous time might provide the best initial guess for the current time level, the complexity of the $k$-means algorithm might overcome the benefit.}

\reminder{To give a complete complexity analysis, we first define the following \textit{on-the-fly-$\alpha$} algorithm, which re-classify all trajectories only every 
$\alpha$ time steps. Therefore, we replace \textbf{Step 2} of the above algorithm with arrival locations given by \{$\mathbf{x}_j(t_{z+i}) | i=1,2,3,...,\alpha$\} for some $z$ a constant multiple of a positive integer $\alpha$ and use the initial centroids 
$$
\frac{1}{\mbox{card}(C_l)} \sum_{\lambda_j=l} \mathbf{X}_j^{z+\alpha}
$$ 
in \textbf{Step 3}.}

\reminder{Assuming that the number of iterations in each $k$-means algorithm is bounded, the overall computational complexity is given by
$$
O\left( \sum_{i=1}^{n/\alpha} M k d \alpha i \right) = O\left( M k d \alpha \sum_{i=1}^{n/\alpha} i \right)
= O\left( \frac{M k d n^2}{\alpha} \right)
$$
where $M$ is the number of trajectories, $k$ is the number of clusters, $d$ is the physical dimensions, and $d\alpha i$ gives the length of each trajectory data at the $i$-th re-classification stage. Therefore, the computational complexity of the \textit{on-the-fly-1} algorithm will be relatively high, given by $O(Mkdn^2)$. Instead, if we know $n$ ahead of the computations, we can determine an $\alpha=O(n)$ given by a certain fraction of $n$. This choice of $\alpha$ leads to a $O(Mkdn)$ algorithm, which has the same complexity as one single $k$-means iteration.}

Like in the previous section, in \textbf{Step 3}, we do not necessarily recompute the mean vector by considering the whole trajectory but may only take the average of the latest position data within each cluster. This improvement will further improve computational efficiency. Also, we might not need to obtain the converged result in any of the intermediate stages if we only care about the clustering result in the last step.

\section{Numerical Examples}
\label{Sec:Examples}
In this section, we will demonstrate our approach in two examples. The first example is the double gyre flow which is a bounded field, and the second velocity field is the van der Pol Duffing oscillator which is an unbounded field. Unless specified otherwise, the WCVE defined below uses SD (\ref{Eqn:WCVE}) rather than MAD (\ref{Eqn:MADE}). 

%%%%%%%%%%%%%%%%%%%%%%%%%%%%%%%%%%%%%%%%%%%%
\subsection{The Double Gyre Flow}
\label{SubSec:double_gyre}

\begin{figure}[!htb]
\includegraphics[width=0.6\textwidth]{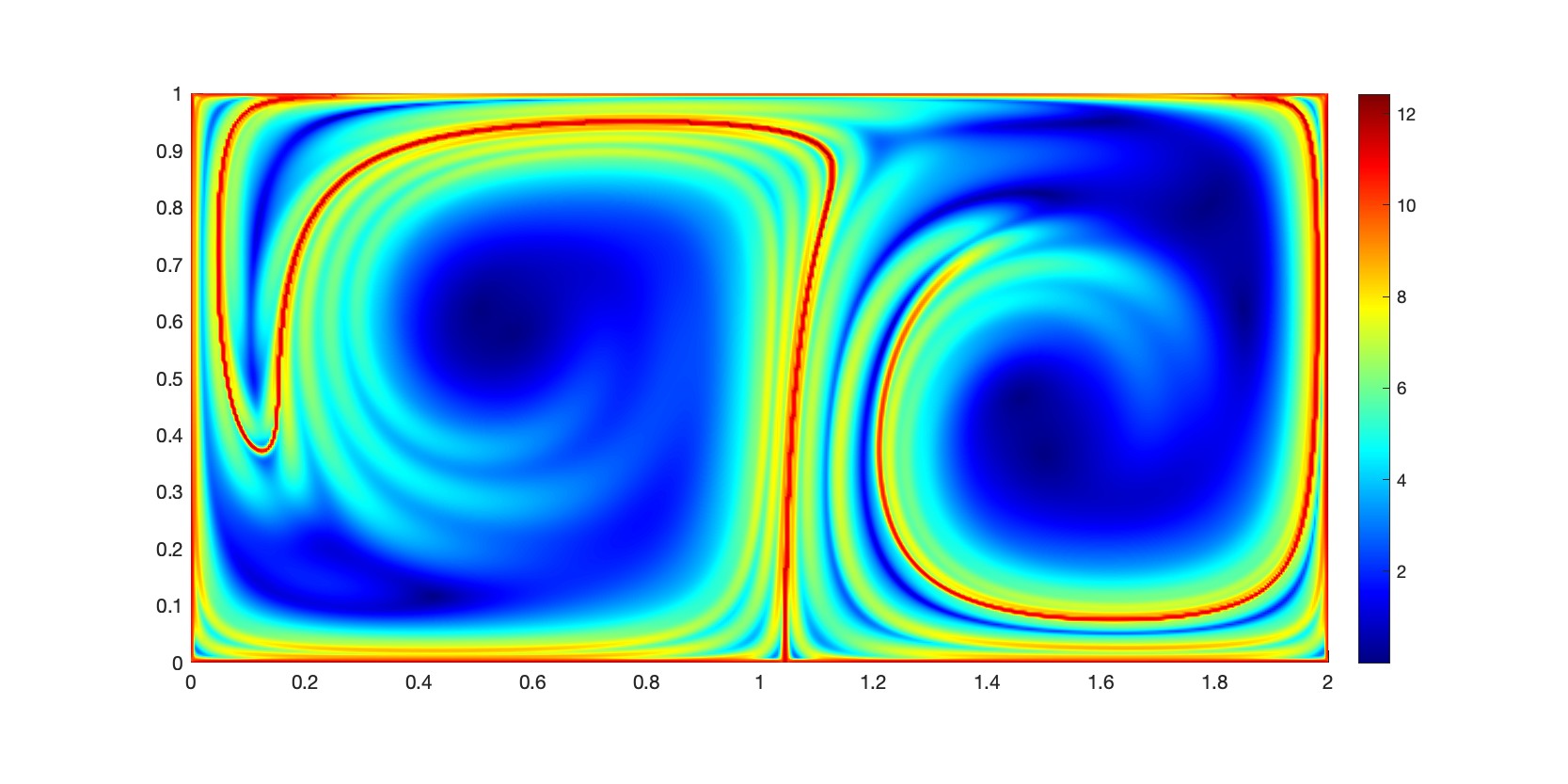}
%(b)\includegraphics[width=0.45\textwidth]{figures/ftlerk5.png}
\caption{\reminder{(Section \ref{SubSec:double_gyre}) The FTLE with $\triangle x = \triangle y = 1/512$ at $T = 15.0$ in  \cite{leu11}.}}
\label{Ex:DoubleGyre1}
\end{figure}

\begin{figure}[!htb]
%(a)\includegraphics[trim=130 50 100 40, clip, width=0.45\textwidth]{figures/tr150b.jpg}
%(b)\includegraphics[trim=130 50 100 40, clip, width=0.45\textwidth]{figures/tr300b.jpg}
%(c)\includegraphics[trim=130 50 100 40, clip, width=0.45\textwidth]{figures/tr450b.jpg}
%(d)\includegraphics[trim=130 50 100 40, clip, width=0.45\textwidth]{figures/tr600b.jpg}
\includegraphics[width=0.95\textwidth]{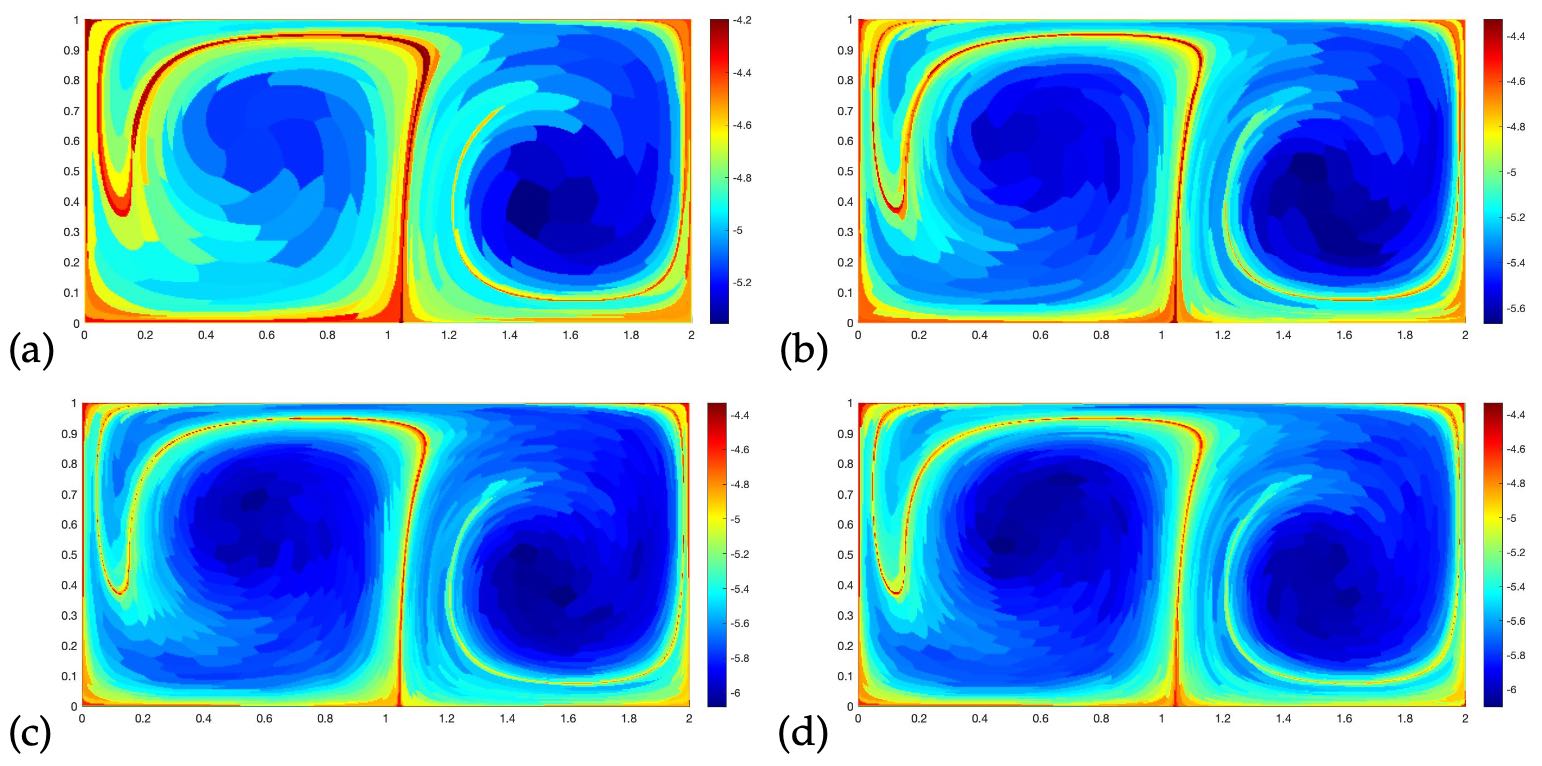}
\caption{(Section \ref{SubSec:double_gyre}) The WCVE computed by  $k$-means algorithm with mesh size $\Delta x = \Delta y$ = 1/256 and $\Delta t$ = 0.1 for $T=15$ with (a) $k=$150, (b) $k=$300, (c) $k=$450 and (d) $k=$600.}
\label{Ex:DoubleGyre_kmean}
\end{figure}

\begin{figure}[!htb]
(a)\includegraphics[trim=130 50 100 40, clip, width=0.45\textwidth]{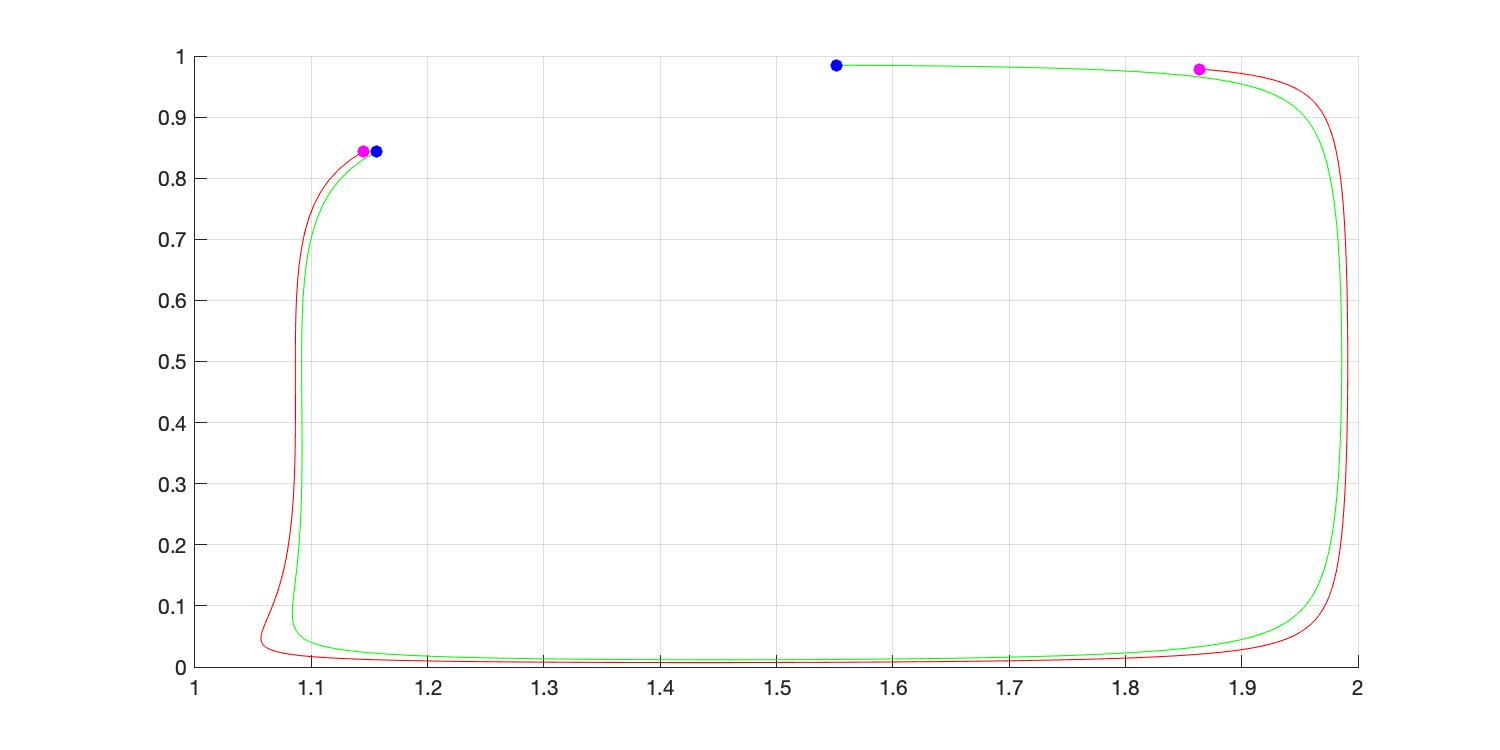}
(b)\includegraphics[trim=130 50 100 40, clip, width=0.45\textwidth]{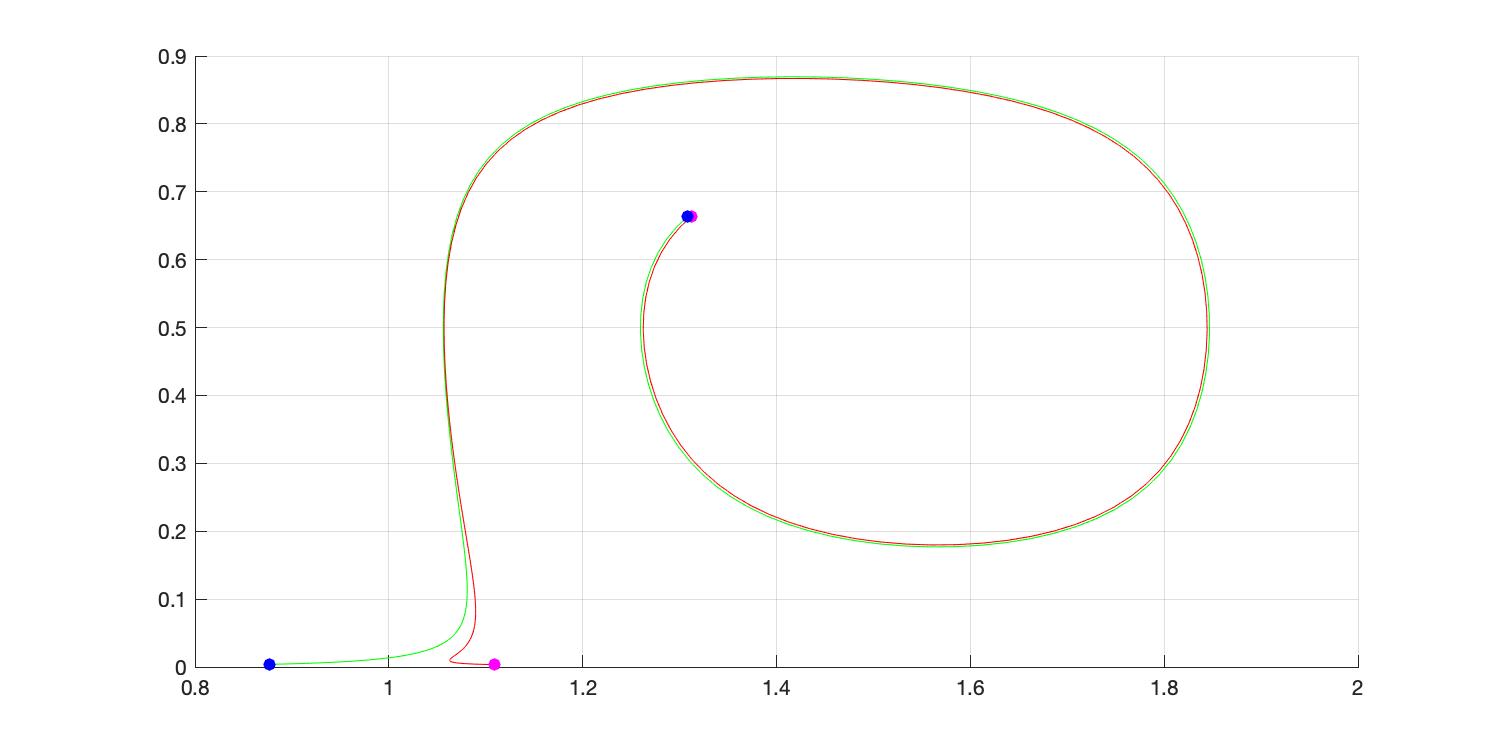}
\caption{(Section \ref{SubSec:double_gyre}) Trajectory of several adjacent particles from (a) (1.15625,0.84375) and (1.14453,0.84375), and (b) (1.30859,0.664062) and (1.3125,0.664062). We see that the FTLE at these locations is large while the corresponding WCVE is small.}
\label{Ex:DoubleGyre_traj}
\end{figure}

\begin{figure}[!htb]
(a)\includegraphics[trim=130 50 100 40, clip, width=0.45\textwidth]{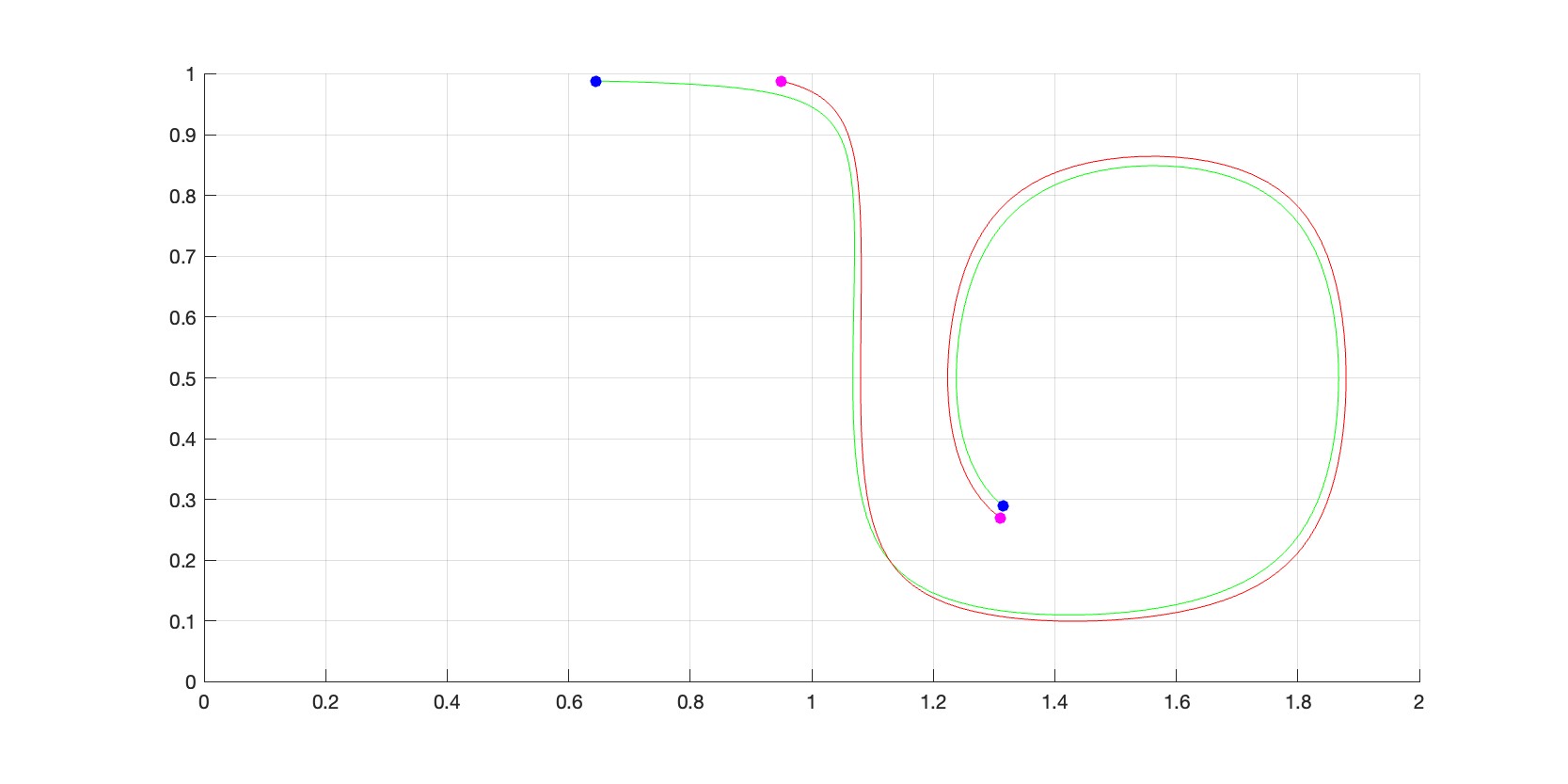}
(b)\includegraphics[trim=130 50 100 40, clip, width=0.45\textwidth]{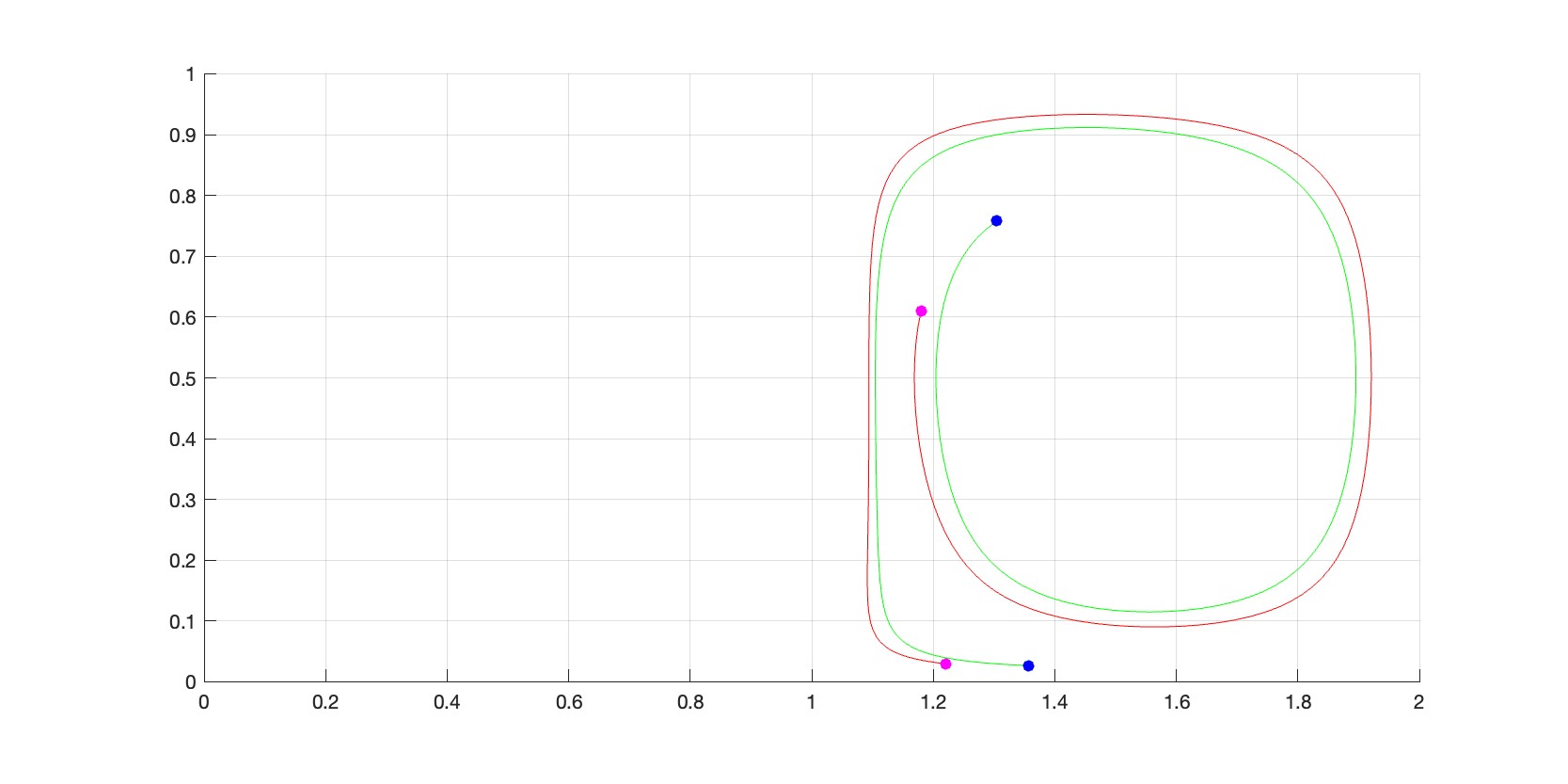}
\caption{\reminder{(Section \ref{SubSec:double_gyre}) Trajectory of several particles from (a) (0.644531,0.988281) and (0.949219,0.988281), and (b) (1.30469,0.757812) and (1.17969,0.609375). Similar history of trajectories with far initial conditions.}}
\label{Ex:DoubleGyre_simtraj}
\end{figure}

\begin{figure}[!htb]
%(a)\includegraphics[trim=130 50 100 40, clip, width=0.45\textwidth]{figures/ran501.jpg}
%(b)\includegraphics[trim=130 50 100 40, clip, width=0.45\textwidth]{figures/ran502.jpg}
%(c)\includegraphics[trim=130 50 100 40, clip, width=0.45\textwidth]{figures/ran503.jpg}
%(d)\includegraphics[trim=130 50 100 40, clip, width=0.45\textwidth]{figures/ran504.jpg}
\includegraphics[width=0.95\textwidth]{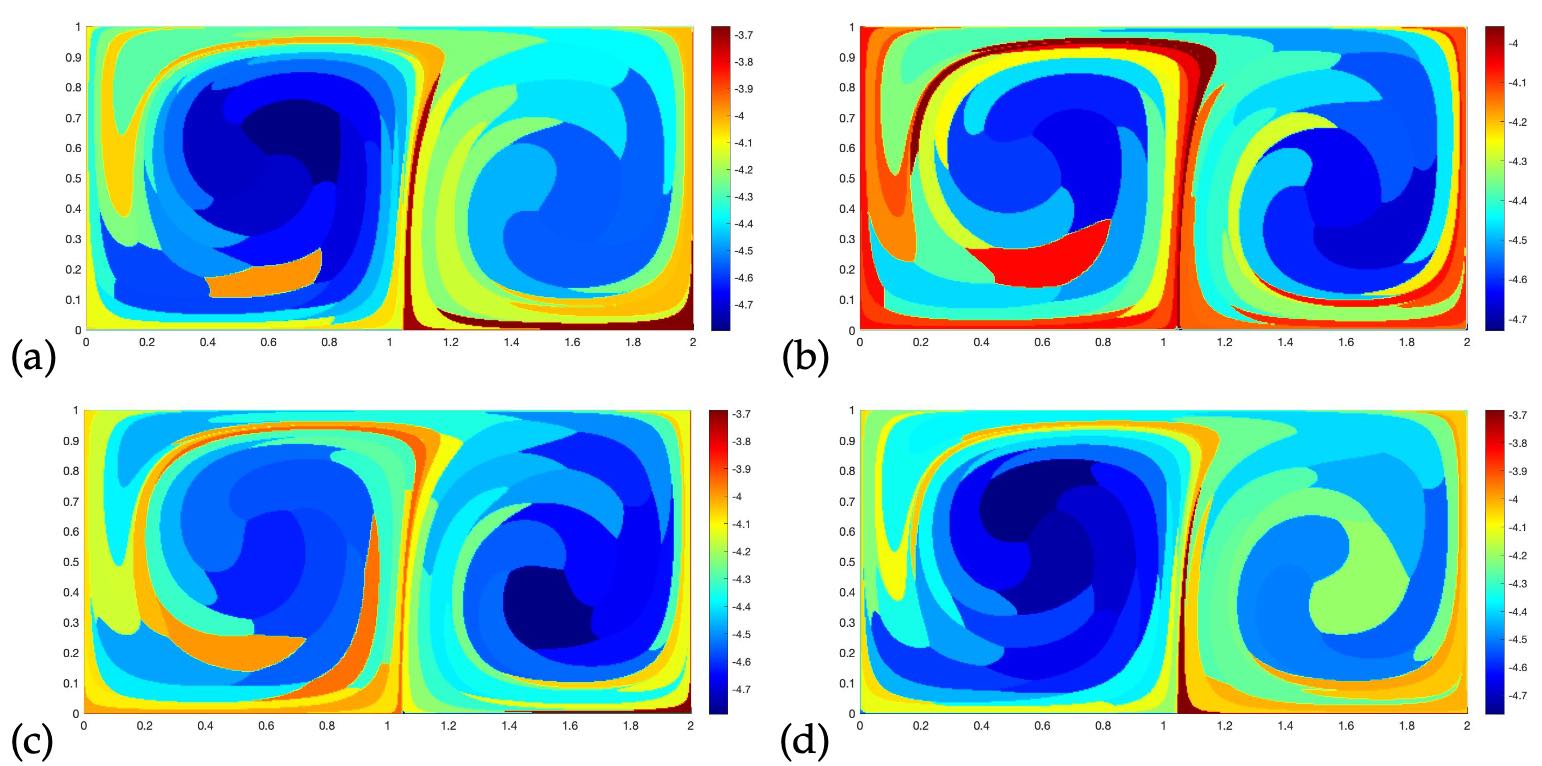}
\caption{(Section \ref{SubSec:double_gyre}) Four independent simulations of WCVE computed by $k$-means algorithm with mesh size $\Delta x = \Delta y$ = 1/256, $\Delta t$ = 0.1 and with $k=50$ at $T=15$.}
\label{Ex:DoubleGyre_ran50}
\end{figure}

\begin{figure}[!htb]
%(a)\includegraphics[trim=130 50 100 40, clip, width=0.45\textwidth]{figures/ran1501.jpg}
%(b)\includegraphics[trim=130 50 100 40, clip, width=0.45\textwidth]{figures/ran1502.jpg}
%(c)\includegraphics[trim=130 50 100 40, clip, width=0.45\textwidth]{figures/ran1503.jpg}
%(d)\includegraphics[trim=130 50 100 40, clip, width=0.45\textwidth]{figures/ran1504.jpg}
\includegraphics[width=0.95\textwidth]{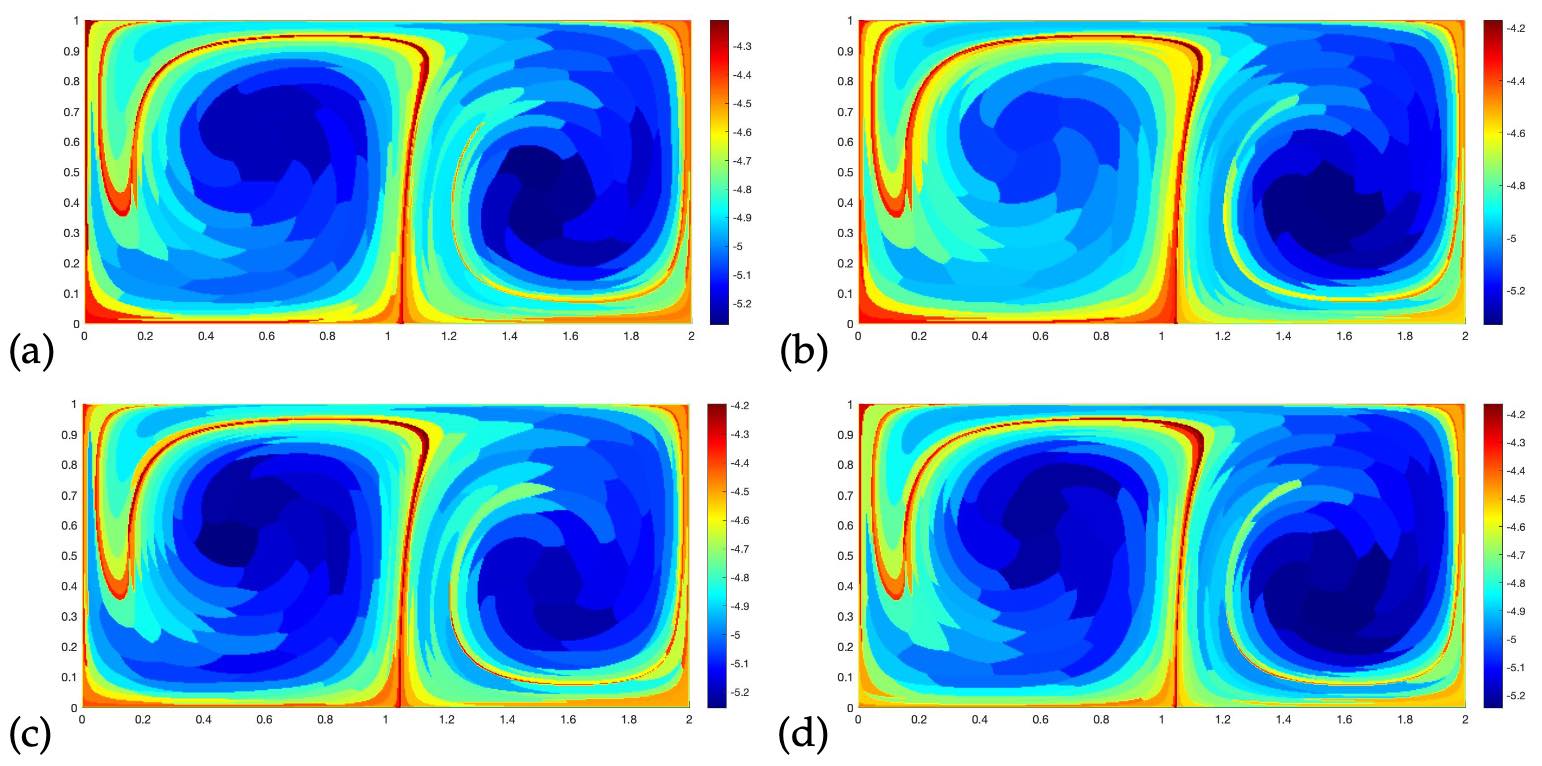}
\caption{(Section \ref{SubSec:double_gyre}) Four independent simulations of WCVE computed by the $k$-means algorithm with mesh size $\Delta x = \Delta y$ = 1/256 , $\Delta t$ = 0.1 and with $k=150$ at $T=15$.}
\label{Ex:DoubleGyre_ran150}
\end{figure}

\begin{figure}[!htb]
%(a)\includegraphics[trim=130 50 100 40, clip, width=0.45\textwidth]{figures/MADE150b.jpg}
%(b)\includegraphics[trim=130 50 100 40, clip, width=0.45\textwidth]{figures/MADE300b.jpg}
%(c)\includegraphics[trim=130 50 100 40, clip, width=0.45\textwidth]{figures/MADE450b.jpg}
%(d)\includegraphics[trim=130 50 100 40, clip, width=0.45\textwidth]{figures/MADE600b.jpg}
\includegraphics[width=0.95\textwidth]{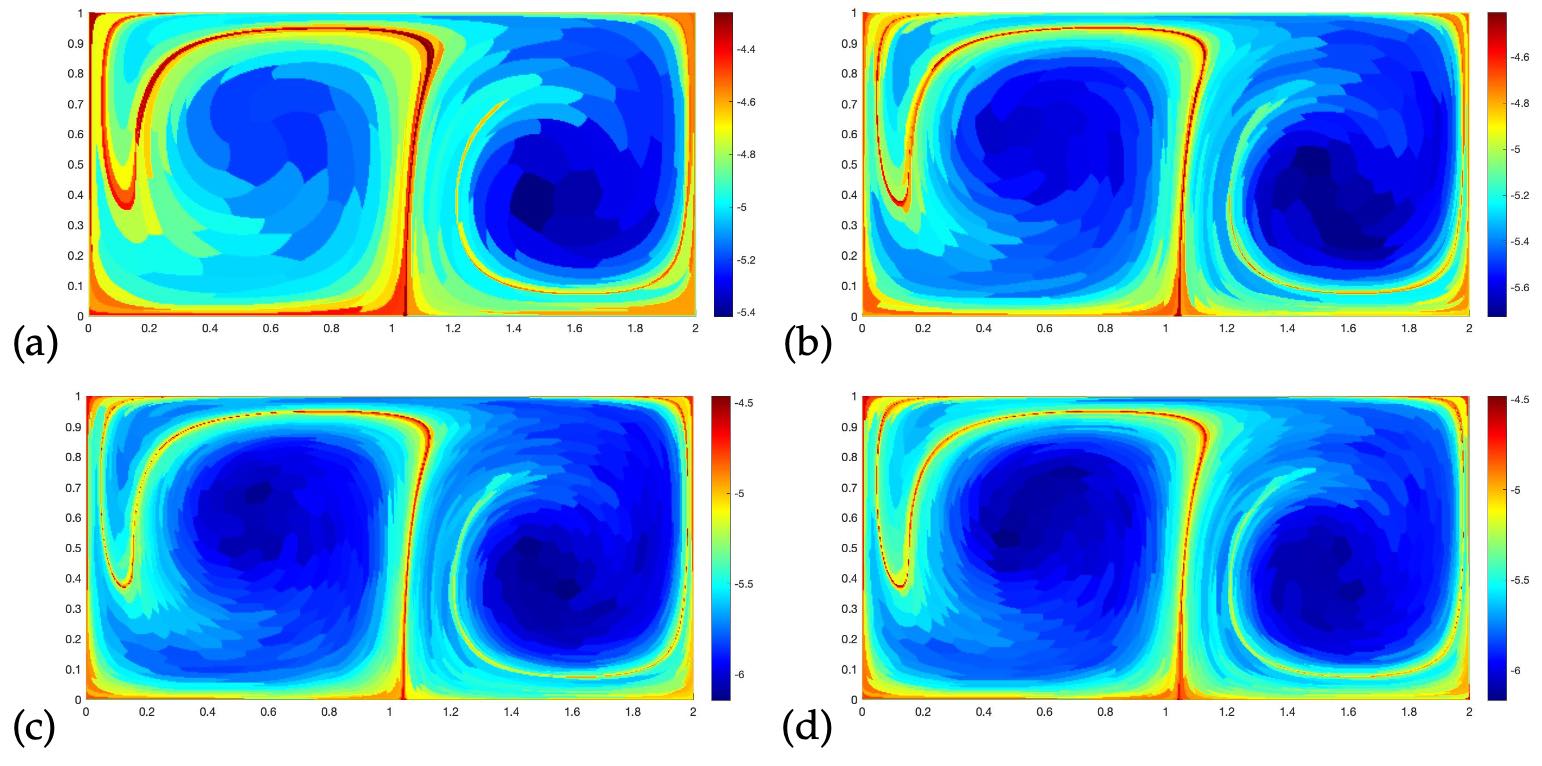}
\caption{(Section \ref{SubSec:double_gyre}) The WCVE defined based on the MAD with mesh size $\Delta x = \Delta y$ = 1/256 and $\Delta t$ = 0.1 for $T=15$ (a) $k=$150, (b) $k=$300, (c) $k=$450 and (d) $k=$600.}
\label{Ex:DoubleGyre_MADE}
\end{figure}

\begin{figure}[!htb]
%(a)\includegraphics[trim=130 50 100 40, clip, width=0.45\textwidth]{figures/fast1150b.jpg}
%(b)\includegraphics[trim=130 50 100 40, clip, width=0.45\textwidth]{figures/fast1300b.jpg}
%(c)\includegraphics[trim=130 50 100 40, clip, width=0.45\textwidth]{figures/fast1450b.jpg}
%(d)\includegraphics[trim=130 50 100 40, clip, width=0.45\textwidth]{figures/fast1600b.jpg}
\includegraphics[width=0.95\textwidth]{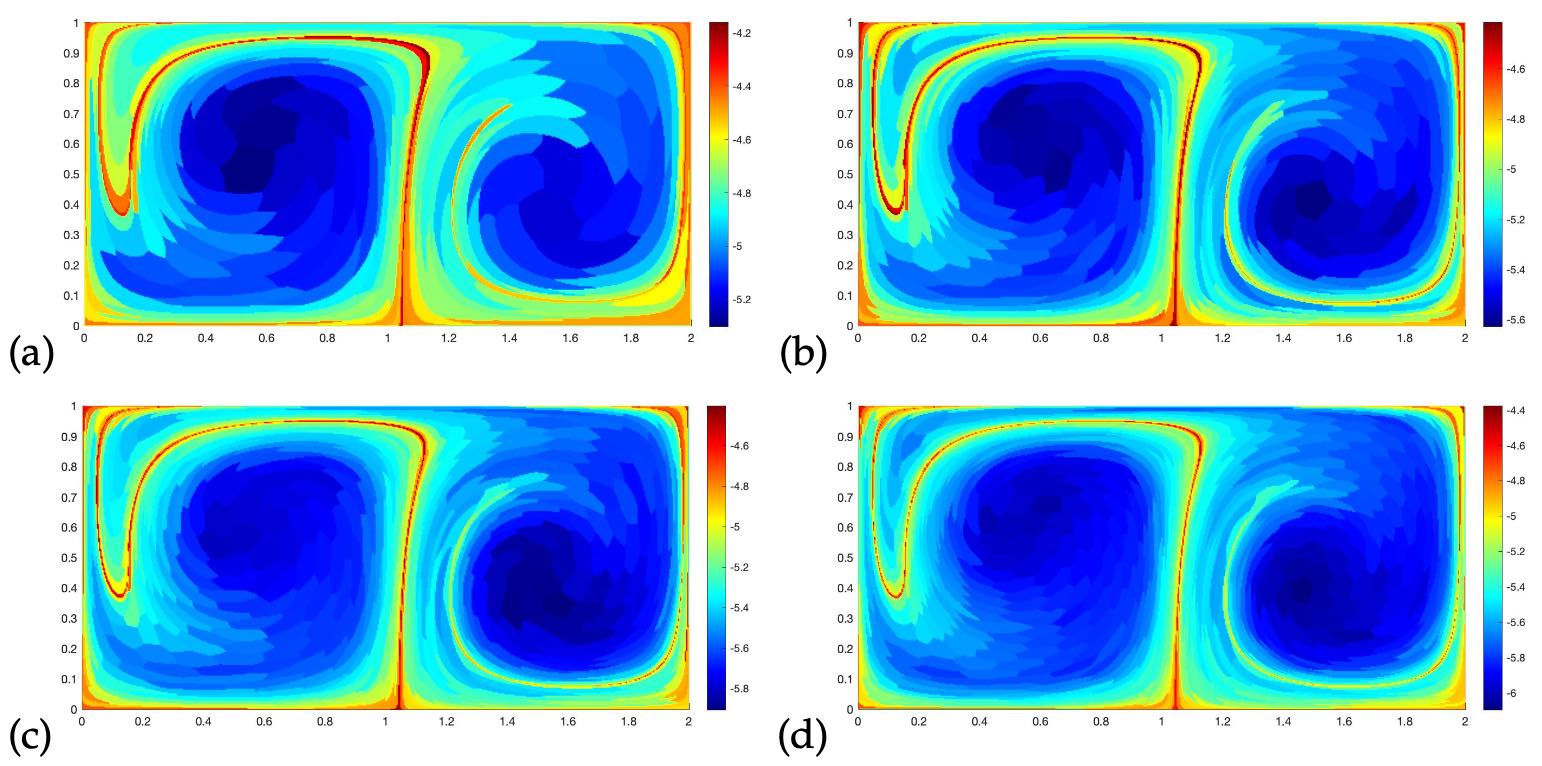}
\caption{(Section \ref{SubSec:double_gyre}) The WCVE computed by an adaptive refinement approach ($N$=4) with mesh size $\Delta x = \Delta y$ = 1/256 and $\Delta t$ = 0.1 for $T=15$ (a) $k=$ 150, (b) $k=$ 300, (c) $k=$ 450 and (d) $k=$ 600.}
\label{Ex:DoubleGyre_adaptive}
\end{figure}

\begin{figure}[!htb]
%(a)\includegraphics[trim=130 50 100 40, clip, width=0.45\textwidth]{figures/300bmid1.jpg}
%(b)\includegraphics[trim=130 50 100 40, clip, width=0.45\textwidth]{figures/300bmid2.jpg}
%(c)\includegraphics[trim=130 50 100 40, clip, width=0.45\textwidth]{figures/300bmid3.jpg}
%(d)\includegraphics[trim=130 50 100 40, clip, width=0.45\textwidth]{figures/300bmid4.jpg}
\includegraphics[width=0.95\textwidth]{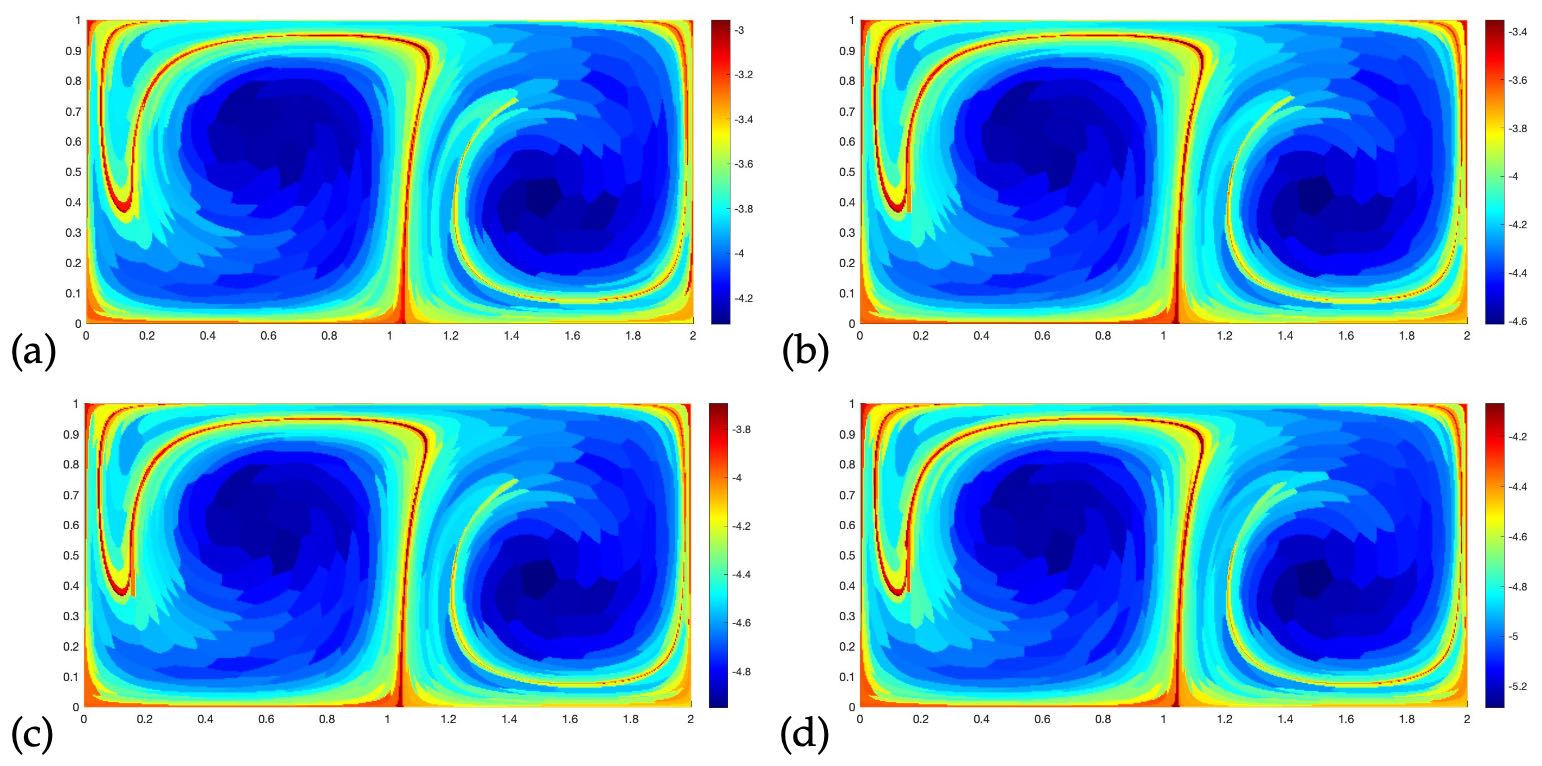}
\caption{(Section \ref{SubSec:double_gyre}) The intermediate WCVE computed by the adaptive refinement approach with $k=300$ for (a) $N=$4, (b) $N=$3, (c) $N=$2 and (d) $N=$1.}
\label{Ex:DoubleGyre_adaptive_inter_300}
\end{figure}

\begin{figure}[!htb]
%(a)\includegraphics[trim=130 50 100 40, clip, width=0.45\textwidth]{figures/450bmid1.jpg}
%(b)\includegraphics[trim=130 50 100 40, clip, width=0.45\textwidth]{figures/450bmid2.jpg}
%(c)\includegraphics[trim=130 50 100 40, clip, width=0.45\textwidth]{figures/450bmid3.jpg}
%(d)\includegraphics[trim=130 50 100 40, clip, width=0.45\textwidth]{figures/450bmid4.jpg}
\includegraphics[width=0.95\textwidth]{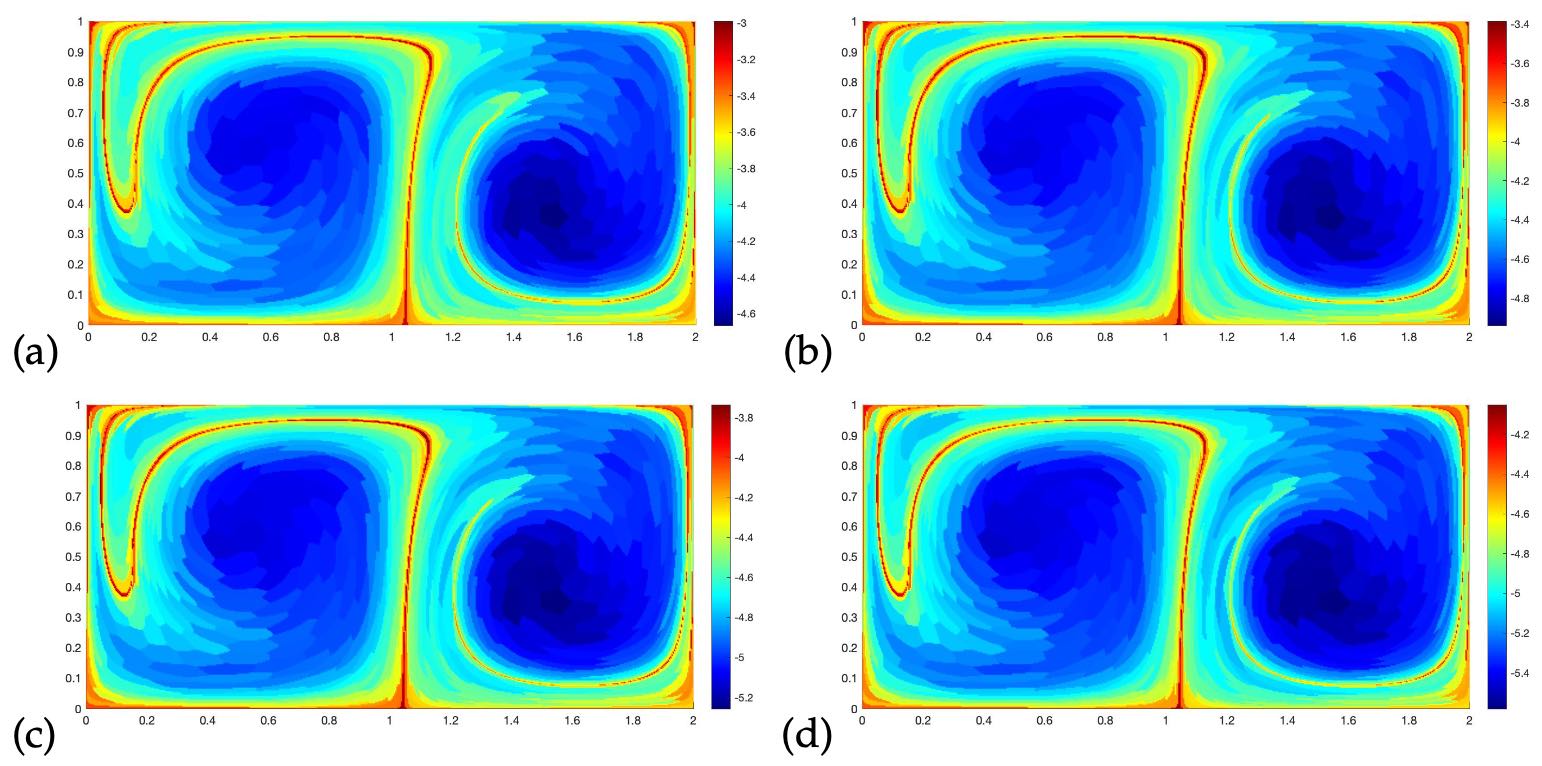}
\caption{(Section \ref{SubSec:double_gyre}) The intermediate WCVE computed by the adaptive refinement approach with $k=450$ for $N=$4, (b) $N=$3, (c) $N=$2 and (d) $N=$1.}
\label{Ex:DoubleGyre_adaptive_inter_450}
\end{figure}

\begin{figure}[!htb]
(a)\includegraphics[trim=130 50 100 40, clip, width=0.45\textwidth]{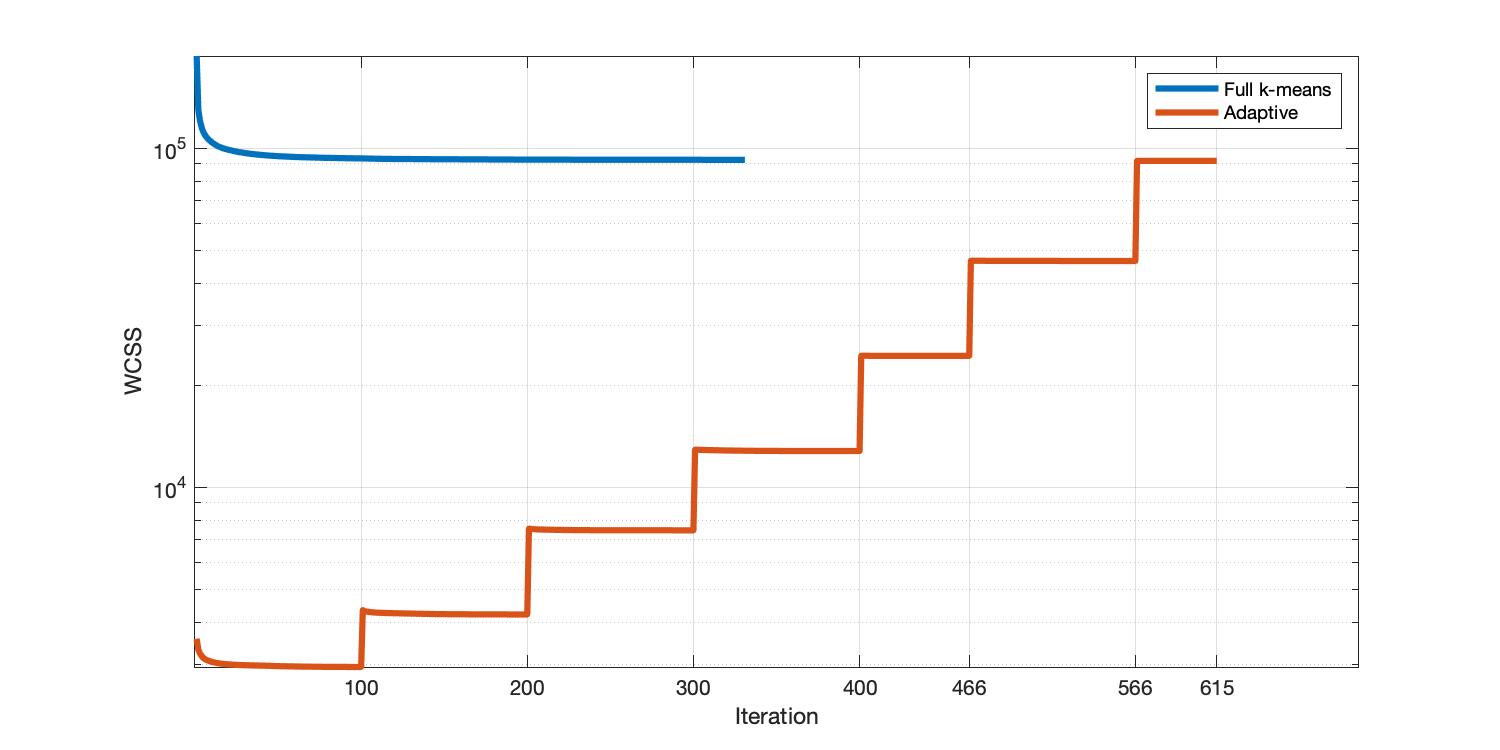}
(b)\includegraphics[trim=130 50 100 40, clip, width=0.45\textwidth]{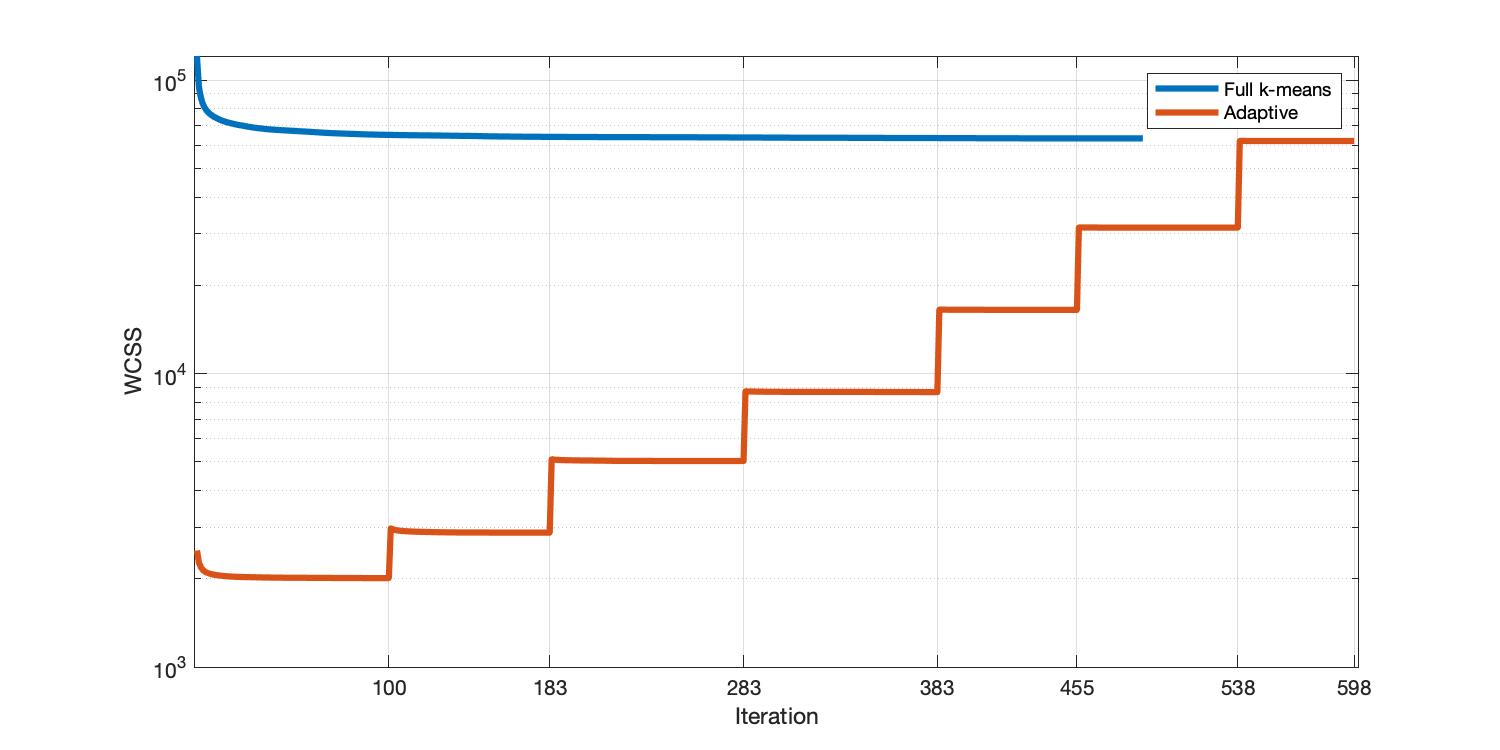}
\caption{(Section \ref{SubSec:double_gyre}) The change in the WCSS (\ref{Eqn:WCSS}) in the full $k$-means iterations and also from the adaptive refinement approach with (a) $k=300$ and (b) $k=450$.}
\label{Ex:DoubleGyre_WCSS}
\end{figure}

\begin{figure}[!htb]
\centering
%(a)\includegraphics[trim=130 70 120 50, clip, width=0.45\textwidth]{figures/larT1500b.jpg}
%(b)\includegraphics[trim=120 220 120 200, clip, width=0.45\textwidth]{figures/LFTLE.jpg}
\includegraphics[width=0.95\textwidth]{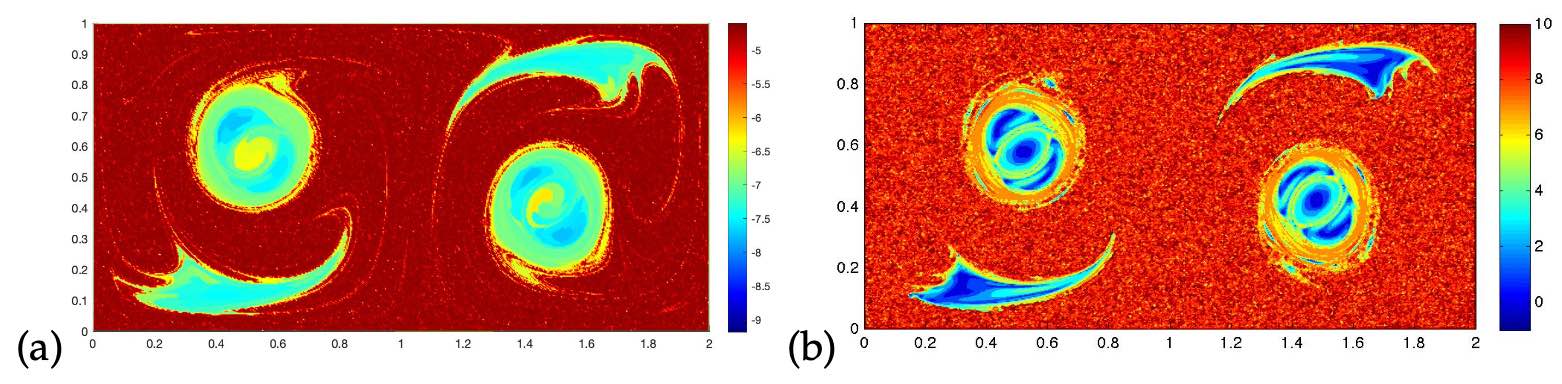}
\caption{(Section \ref{SubSec:double_gyre}) (a) The WCVE field computed using the adaptive refinement approach ($N$=10) with $\Delta x = \Delta y$ = 1/256 and $\Delta t$ = 0.1 for $T$= 320 for $k$=1500. (b) The FTLE field with $\Delta x = \Delta y$ = 1/256 for $T$= 320 taken in \cite{leu13}.}
\label{Ex:DoubleGyre_larT_FTLE}
\end{figure}

\begin{figure}[!htb]
%(a)\includegraphics[trim=130 50 100 40, clip, width=0.45\textwidth]{figures/not_0.jpg}
%(b)\includegraphics[trim=130 50 100 40, clip, width=0.45\textwidth]{figures/not_1.jpg}
\includegraphics[width=0.95\textwidth]{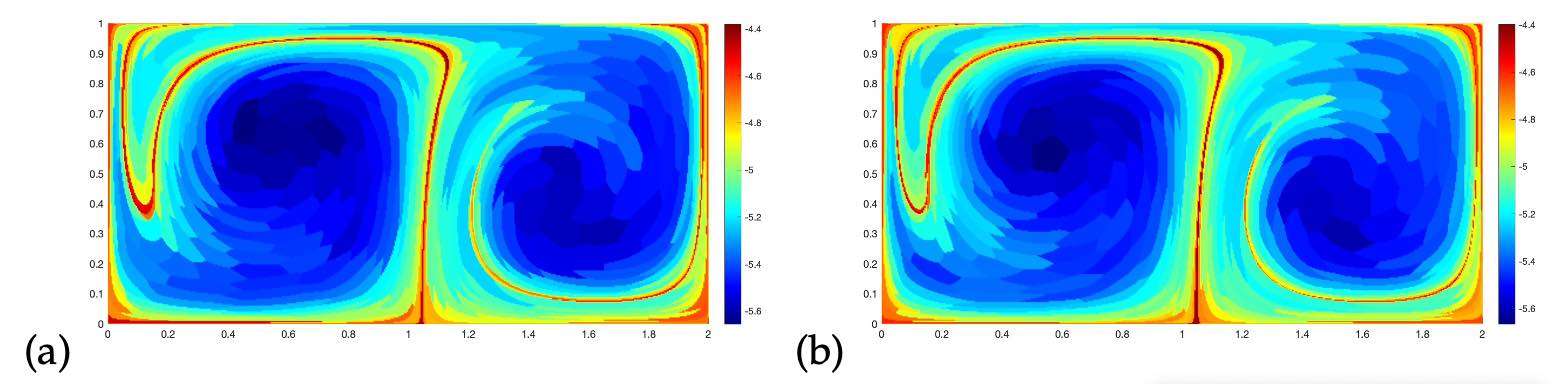}
\caption{\reminder{(Section \ref{SubSec:double_gyre}) The WCVE field with mesh size $\Delta x = \Delta y$ = 1/256 , $\Delta t$ = 0.1 and with $k=300$ at $T=15$ and computed using (a) the full $k$-means algorithm and (b) the adaptive refinement approach with the same initial condition.}}
\label{Ex:DoubleGyre_con}
\end{figure}

\begin{figure}[!htb]
%(a)\includegraphics[trim=130 50 100 40, clip, width=0.45\textwidth]{figures/51303b_ada_inde1.jpg}
%(b)\includegraphics[trim=130 50 100 40, clip, width=0.45\textwidth]{figures/51303b_ada_inde2.jpg}
%(c)\includegraphics[trim=130 50 100 40, clip, width=0.45\textwidth]{figures/51303b_ada_inde3.jpg}
%(d)\includegraphics[trim=130 50 100 40, clip, width=0.45\textwidth]{figures/51303b_ada_inde4.jpg}
\includegraphics[width=0.95\textwidth]{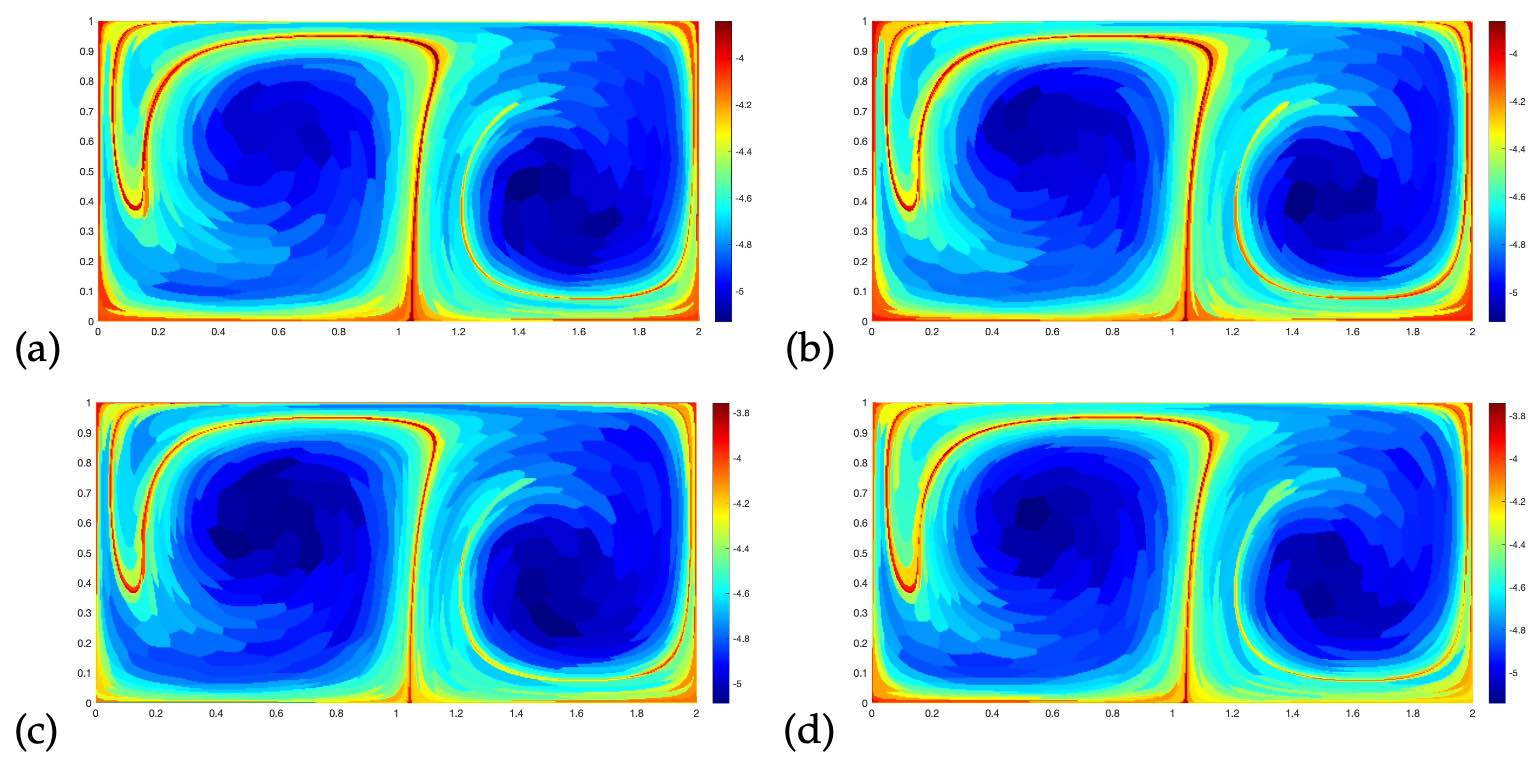}
\caption{\reminder{(Section \ref{SubSec:double_gyre}) Four independent simulations of WCVE computed by the adaptive refinement approach ($N$=5) with mesh size $\Delta x = \Delta y$ = 1/256 , $\Delta t$ = 0.3 and with $k=300$ at $T=15$.}}
\label{Ex:DoubleGyre_stable}
\end{figure}

\begin{figure}[!htb]
%(a)\includegraphics[trim=130 50 100 40, clip, width=0.45\textwidth]{figures/300bt10.jpg}
%(b)\includegraphics[trim=130 50 100 40, clip, width=0.45\textwidth]{figures/300bt102.jpg}
%(c)\includegraphics[trim=130 50 100 40, clip, width=0.45\textwidth]{figures/300bt103.jpg}
%(d)\includegraphics[trim=130 50 100 40, clip, width=0.45\textwidth]{figures/300bt104.jpg}
\includegraphics[width=0.95\textwidth]{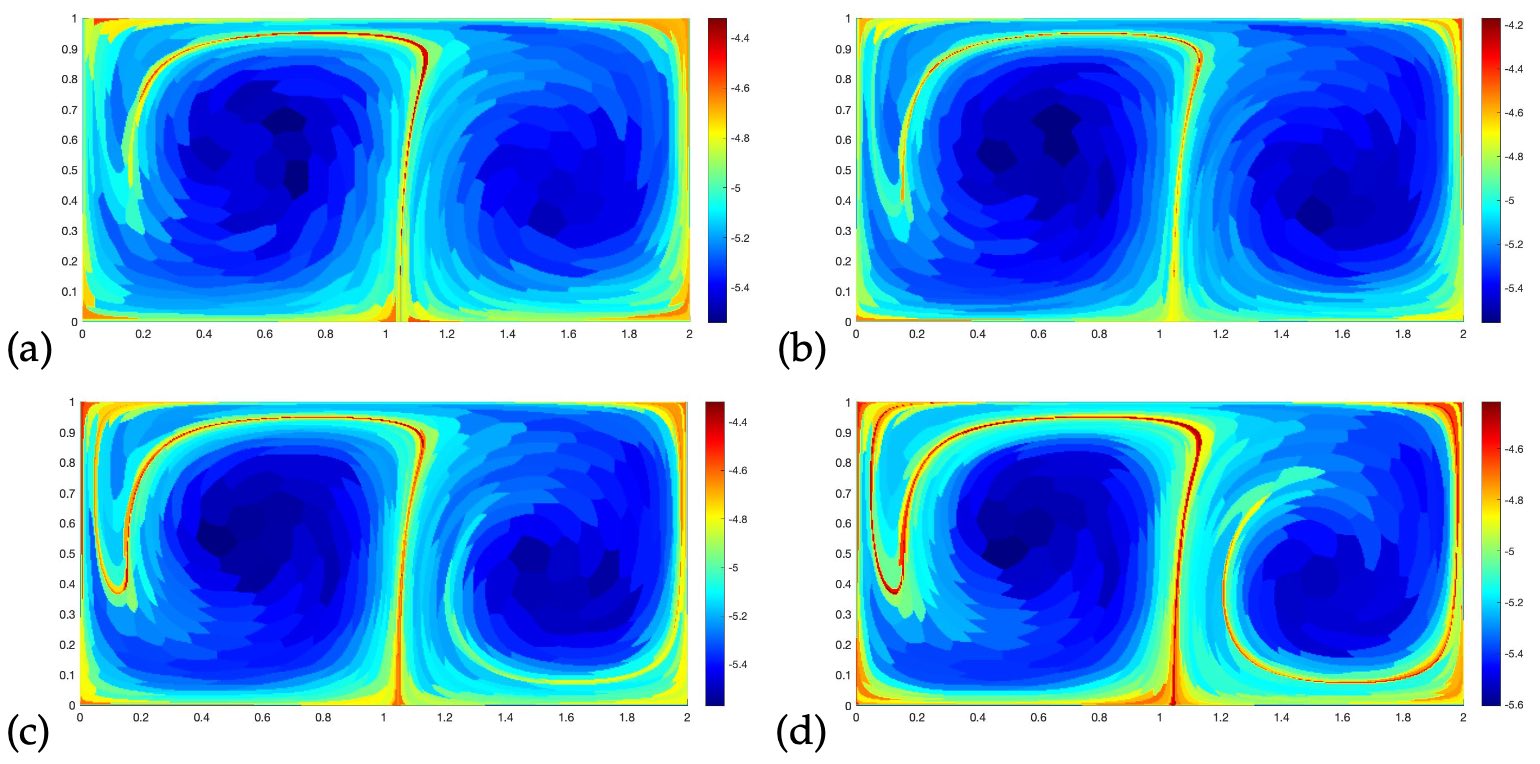}
\caption{(Section \ref{SubSec:double_gyre}) The WCVE computed by \textit{on-the-fly} approach with $k=300$ for (a) $t_z$= 10, (b) $t_z$= 11.5, (c) $t_z$= 13 and (d) $t_z$= 15.}
\label{Ex:DoubleGyre_fly}
\end{figure}

\begin{figure}[!htb]
%(a)\includegraphics[trim=130 50 100 40, clip, width=0.45\textwidth]{figures/onfly5-300.jpg}
%(b)\includegraphics[trim=130 50 100 40, clip, width=0.45\textwidth]{figures/onfly10-300.jpg}
\includegraphics[width=0.95\textwidth]{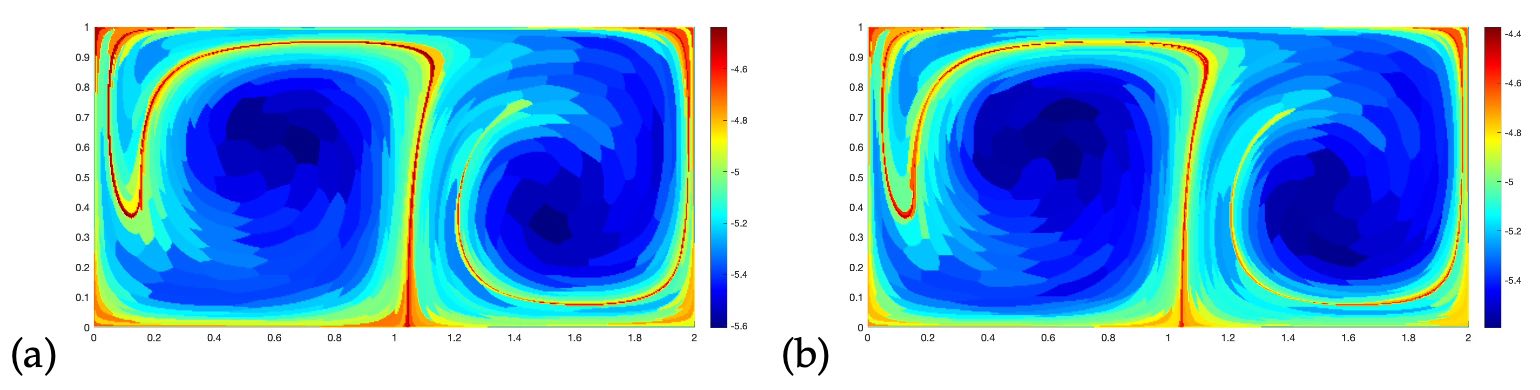}
\caption{\reminder{(Section \ref{SubSec:double_gyre}) The WCVE with $k=300$ computed by (a) \textit{on-the-fly-5}, (b) \textit{on-the-fly-10}.}}
\label{Ex:DoubleGyre_fly5_10}
\end{figure}

\begin{table}[!htb]
\centering
(a)
\begin{tabular}{|c|c|c|c|c|c|c|}
\hline
$k$& 150 & 300 & 450 & 600 \\ 
\hline\hline
Full $k$-means & 960 & 1961 & 3008 & 3115 \\ 
Adaptive & 337 & 481 & 609 & 607 \\ 
\reminder{\textit{on-the-fly-5}} & 749 & 1207 & 1590 & 1914\\
\reminder{\textit{on-the-fly-10}} & 710 & 945 & 1263 & 1481\\
\hline 
\end{tabular} \\
\vspace{0.3cm}
(b)
\begin{tabular}{|c|c|c|c|c|c|c|} 
\hline
$k$& 150 & 300 & 450 & 600 \\ 
\hline\hline
Full $k$-means  & 178219 & 93483 & 64150 & 48998 \\ 
Adaptive & 177416 & 92198 & 61991 & 46502 \\ 
\reminder{\textit{on-the-fly-5}} & 177667 & 92140 & 62097 & 46651\\
\reminder{\textit{on-the-fly-10}} & 177292 & 92339 & 62105 & 46683\\
\hline 
\end{tabular}
\caption{(Section \ref{SubSec:double_gyre}) Comparisons of the full $k$-means algorithm and the adaptive refinement algorithm with measurements taken from the average of 50 trials. (a) The average CPU time (seconds) and (b) the averaged minimized WCSS (\ref{Eqn:WCSS}). We have rounded all numbers to the closest integer value.}
\label{Ex:DoubleGyre_TimeError}
\end{table}

This example is taken from \cite{shalekmar05} to describe a periodically varying double-gyre. The flow is modeled by the following stream-function $\psi(x,y,t)=A \sin[ \pi g(x,t) ] \sin(\pi y)$ where
\begin{eqnarray*}
g(x,t) &=& a(t) x^2 + b(t) x \, , \nonumber\\
a(t) &=& \epsilon \sin(\omega t) \, , \nonumber\\
b(t) &=& 1- 2\epsilon \sin(\omega t) \, .
\end{eqnarray*}
In this example, we follow \cite{shalekmar05,leu11} and use $A=0.1$, $\omega=2\pi/10$ and $\epsilon=0.1$. 
In all computations, we partition the domain [0,2]$\times$[0,1] with $\Delta x = \Delta y$ = 1/256 to compute the WCVE. As a reference, we have shown in Figure \ref{Ex:DoubleGyre1} the FTLE solution computed using the approach as discussed in \cite{leu11}.

\paragraph{WCVE with the Full $k$-means Algorithm.} Figure \ref{Ex:DoubleGyre_kmean} shows the WCVE solution computed using the full $k$-means clustering algorithm with different number of clusters $k$ ranging from 150 to 600. The main observation is that the structure in the WCVE shows some similarities with that of the FTLE. Both quantities show a distinguish filament along $x=1$, which stretches towards the left half of the domain. The solution from the WCVE has fewer features. Because the computation of WCVE relies on the overall behavior of the particle trajectory, it seems that the relatively large value of the FTLE at those locations might come from the variation of the particle motion near the end of the time period. To further investigate this observation, we have shown the trajectories of several particles in Figure \ref{Ex:DoubleGyre_traj}. In these subplots, we have identified two particles that are closed at the initial time, have been clustered into the same group in the WCVE calculations, and the FTLE values are significant. In both cases, the separation between two particles is mostly minor but becomes significant only for the final time interval. In this case, the FTLE is large, but WCVE is relatively insufficient. Therefore, our proposed WCVE depends on a more global-in-time property of the particle trajectories, while the FTLE depends solely on the change in the separation at the final time.

The second feature in the WCVE is that the filament structure is getting sharper and thinner as we increase the number of clusters in the $k$-means algorithm. Because trajectories near the large WCVE region are generally less similar, we can put these data into different groups as we increase $k$, which effectively lower the WCVE values near these regions. We can also see that the most significant values of the WCVE are reduced too when comparing Figure \ref{Ex:DoubleGyre_kmean}(a) and (d).

\reminder{The third property we would like to mention is that our proposed WCVE can also measure the similarity in trajectories from nonlocal initial conditions. Even though the initial location could be quite different, our approach might still classify these particles into the same group if they tend to move toward each other later. In Figure \ref{Ex:DoubleGyre_simtraj}, we have shown two pairs of trajectories starting at different regions in the domain but have been classified into the same cluster by our algorithms with $k=450$.}

\paragraph{The Effect on the $k$-means Initialization.} Next, we would like to investigate the effect of the initialization in the $k$-means algorithm on the WCVE computations. Figure \ref{Ex:DoubleGyre_ran50} and Figure \ref{Ex:DoubleGyre_ran150} show the computed WCVE field with four independent trials with $k=50$ and $k=150$, respectively. Since the current implementation of $k$-means clustering uses a randomized initial guess, the difference in these solutions implies that our WCVE depends on the initial guess of the $k$-means. However, as we increase the number of clusters, the computed solutions seem to be less dependent on the initial guess of the clustering algorithm, and the change in the WCVE structure is less significant. Therefore, a large $k$ seems to ensure the localization and the stability of the WCVE field.

\paragraph{The WCVE based on SD and MAD.} In Figure \ref{Ex:DoubleGyre_MADE}, we computed the WCVE using the full $k$-means algorithm with the definition replaced by MAD (\ref{Eqn:MADE}). Comparing these results with those in Figure \ref{Ex:DoubleGyre_kmean}, we do not see any significant difference in the solutions based on these two different measures. They generally do not show any distinguishable features, but both can identify similar coherent structures in the dynamical system.

\paragraph{The Adaptive Refinement Algorithm.} Figure \ref{Ex:DoubleGyre_adaptive} shows the solutions computed using the adaptive refinement algorithm with the same setting as in Figure \ref{Ex:DoubleGyre_kmean}. We observe that these solutions are indistinguishable from each other. To have some more quantitative comparison, we also look into the performance of the two approaches. We have recorded the CPU time and the WCSS (\ref{Eqn:WCSS}) for the $k$-means algorithm and adaptive refinement approach. All simulations in this work are computed using a laptop computer with a 2.4 GHz Intel core $i5$ processor. In Table \ref{Ex:DoubleGyre_TimeError}(a), we see that the CPU time for the adaptive refinement approach is significantly shorter (to approximately one-fifth) than that of the full $k$-means algorithm. The clustering performance has also been improved since the minimized WCSS for the adaptive refinement approach is less than the $k$-means algorithm, as shown in Table \ref{Ex:DoubleGyre_TimeError}(b).

Figure \ref{Ex:DoubleGyre_adaptive_inter_300} and Figure \ref{Ex:DoubleGyre_adaptive_inter_450} show the intermediate solutions of WCVE computed by the adaptive refinement approach, with $k=300$ and $k=450$. Comparing subplots (a) and (d), we see that the coarsest level $N=4$ in (a) already provides some rough structures of the WCVE. As we provide more information in-between different levels, these structures become more sharpened and more obvious. In general, if all flow trajectories are already well-captured by a given time discretization, the further refinement in the time sampling will not provide more information to tell trajectories apart. We have a similar observation in the plot of the change in the WCSS, as shown in Figure \ref{Ex:DoubleGyre_WCSS} where we have collected the WCSS data using $k=300$ and $450$ from the full $k$-means algorithm and also our proposed adaptive refinement approach starting with $N=6$. For the full $k$-means algorithm (plotted in solid blue lines), we see that the WCSS in both cases converges nicely. For the adaptive refinement case, the sudden jumps in the WCSS come from the upsampling of the data from low to high dimensions. Even though the WCSS after each jump further reduces as iterations go, the change is rather minimal, and the curve seems flat until the next upsampling step.

Because of the adaptive algorithm, we can obtain a long-time computation efficiently. In this example, we increase the time period from $T=10$ to $T=320$, equivalent to $32$ periods of the oscillating perturbation. In Figure \ref{Ex:DoubleGyre_larT_FTLE}(a), we choose $N = 10$ and a significant value of $k=1500$ (still much smaller than $M=257\times513=131841$) to cluster the complicated trajectories. As a comparison, we show in Figure \ref{Ex:DoubleGyre_larT_FTLE}(b) the FTLE computed using the backward phase flow method as discussed in \cite{leu13}. Similar to the FTLE solution, the value is relatively large and homogeneous in most regions. This observation indicates that trajectories are complicated, and any small perturbation in the initial location will grow significantly in time. We can also identify similar regions in the domain where the within-cluster variability is small. Particles from those regions seem to travel in a much similar pattern. Because of the adaptive algorithm, we see a great improvement in computational time. The full $k$-means algorithm took 51760 seconds, while the adaptive algorithm with $N=10$ required only approximately one-third of the computational time (16462 seconds) to complete.

\reminder{Before we end this discussion, we consider both the stability and robustness adaptive algorithm. Although Table \ref{Ex:DoubleGyre_TimeError}(b) already showed that the minimized WCSS from the full $k$-means algorithm is generally different from that by the adaptive algorithm, we want to investigate how different the classification results would be. Figure \ref{Ex:DoubleGyre_con} shows the solution by these two algorithms with the same initial centroids in the $k$-means clustering with $k=300$. We see that these solutions match qualitatively well and do not significantly differ in the computed WCVE structure. We are also interested in the effect of $\Delta t$ on the classification results. When $\Delta t$ is small enough to ensure the stability of the numerical scheme for solving the corresponding ODE, we see that the WCVE solution is robust and does not show any significant difference as we change the initial centroids in the $k$-means algorithm. Figure \ref{Ex:DoubleGyre_stable} show the four independent solutions of WCVE computed by the adaptive refinement approach with $k=300$ and relatively large time steps $\Delta t=0.3$. Once again, these solutions match qualitatively well and do not significantly differ in the computed WCVE structure.}

\paragraph{The \textit{On-the-fly} Algorithm.}
Finally, we consider the \textit{on-the-fly} algorithm in Figure \ref{Ex:DoubleGyre_fly}. The figure shows the WCVE at different times. Instead of re-starting the whole computations, the algorithm updates the solution when we have more measurements. We see that the solution at the final time $T=15$ well matches with those in Figure \ref{Ex:DoubleGyre_kmean}(b) or Figure \ref{Ex:DoubleGyre_adaptive}(b). \reminder{Figure \ref{Ex:DoubleGyre_fly5_10} shows the solutions computed with \textit{on-the-fly-5} and \textit{on-the-fly-10}. We see that these solutions also qualitatively match well with Figure \ref{Ex:DoubleGyre_fly}(d) computed using \textit{on-the-fly-1}.} \reminder{Table \ref{Ex:DoubleGyre_TimeError} also shows some quantitative measurements of these \textit{on-the-fly} algorithms. In Table \ref{Ex:DoubleGyre_TimeError}(b), we see that the minimized WCSS from both \textit{on-the-fly} algorithms are different from those by the full $k$-means and the adaptive algorithm, which indicates that the final classified solutions are indeed different. Table \ref{Ex:DoubleGyre_TimeError}(a) shows the computational time of our \textit{on-the-fly} algorithm. We see that this algorithm can reduce the convergence in the $k$-means clustering because we can always provide reasonably good initial centroids for the classification step.}

%%%%%%%%%%%%%%%%%%%%%%%%%%%%%%%%%%%%%%%%%%%%
\subsection{The Forced-damped Duffing van del Pol Equation}
\label{SubSec:Duffing van del Pol}

\begin{figure}[!htb]
\centering
\includegraphics[trim=10 20 40 10, clip, width=0.45\textwidth]{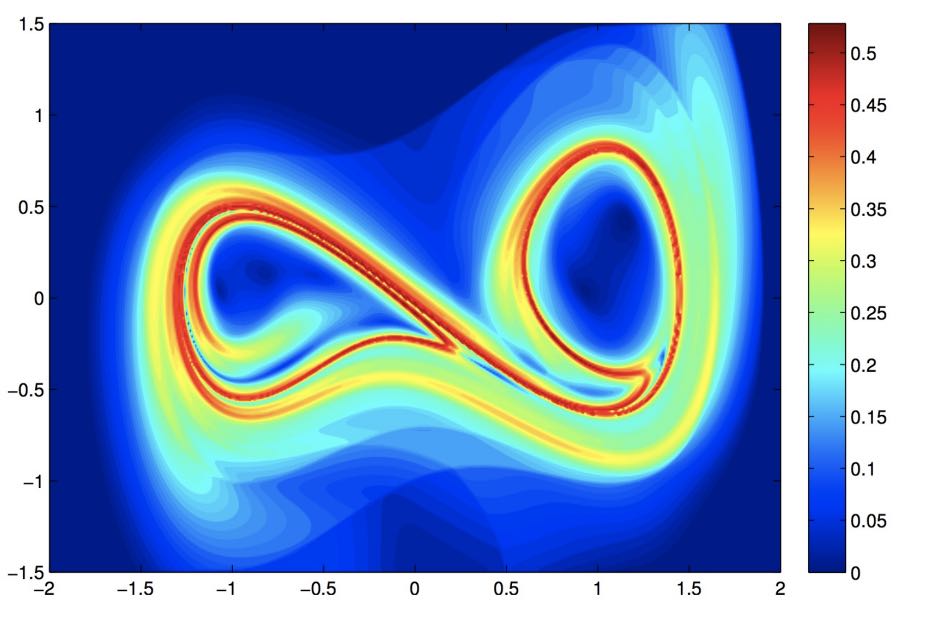}
\caption{(Section \ref{SubSec:Duffing van del Pol}) The FTLE field computed using the proposed Eulerian approach with $\Delta x = \Delta y$ = 0.01 in \cite{youwonleu17}. The negative FTLE values are replaced by 0.}
\label{Ex:Duffing1}
\end{figure}

\begin{figure}[!htb]
%(a)\includegraphics[trim=250 50 250 40, clip, width=0.45\textwidth]{figures/tr100u.jpg}
%(b)\includegraphics[trim=250 50 250 40, clip, width=0.45\textwidth]{figures/tr250u.jpg}
%(c)\includegraphics[trim=250 50 250 40, clip, width=0.45\textwidth]{figures/tr400u.jpg}
%(d)\includegraphics[trim=250 50 250 40, clip, width=0.45\textwidth]{figures/tr550u.jpg}
\includegraphics[width=0.95\textwidth]{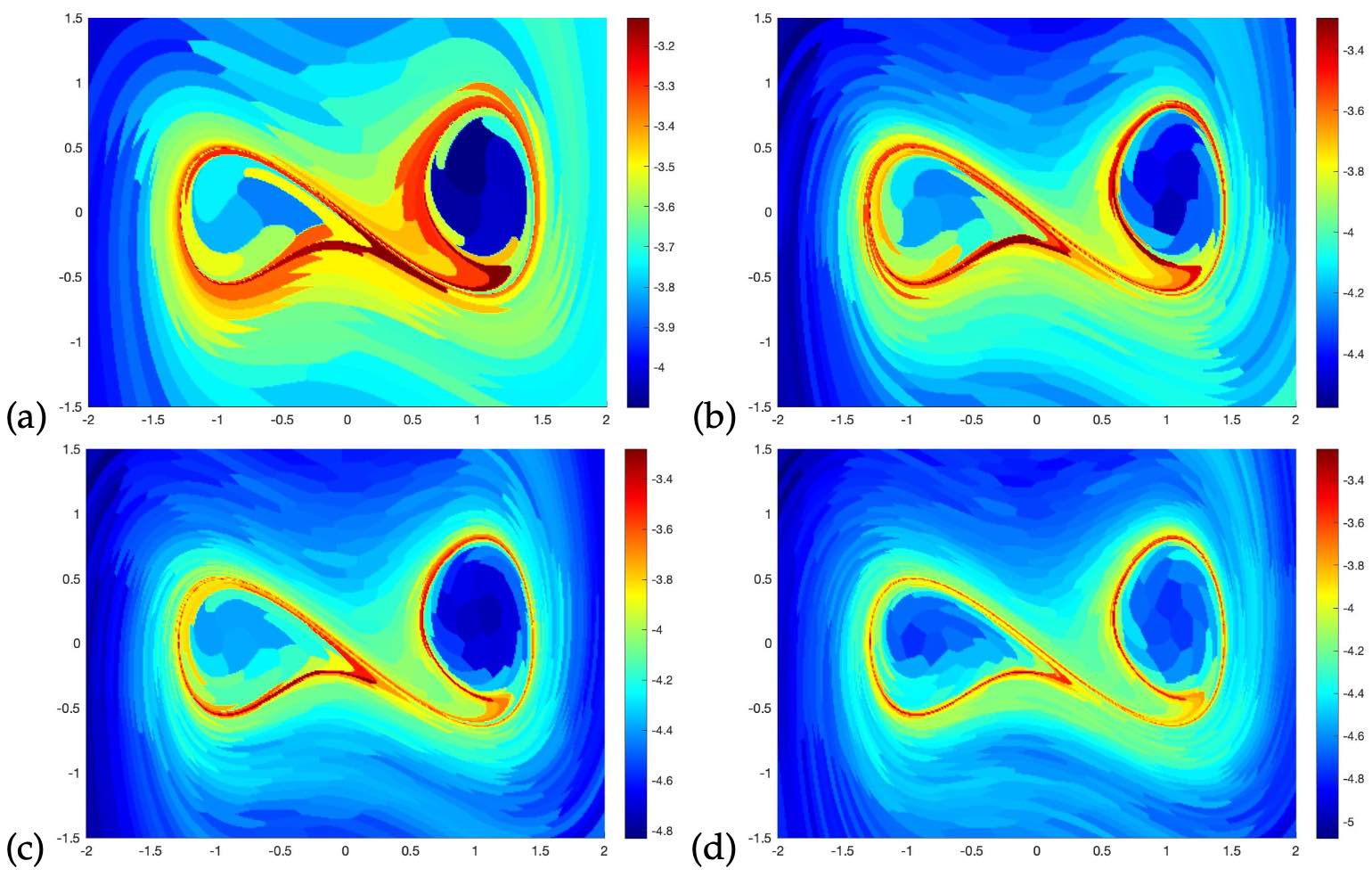}
\caption{(Section \ref{SubSec:Duffing van del Pol}) The WCVE computed by $k$-means algorithm with mesh size $\Delta x = \Delta y$ = 1/100 and $\Delta t$ = 0.1 for $T$=10 with (a) $k=$ 100, (b) $k=$ 250, (c) $k=$ 400, and (d) $k=$ 550.}
\label{Ex:Duffing_kmean}
\end{figure}

\begin{figure}[!htb]
%(a)\includegraphics[trim=250 50 250 40, clip, width=0.45\textwidth]{figures/fast1100u.jpg}
%(b)\includegraphics[trim=250 50 250 40, clip, width=0.45\textwidth]{figures/fast1250u.jpg}
%(c)\includegraphics[trim=250 50 250 40, clip, width=0.45\textwidth]{figures/fast1400u.jpg}
%(d)\includegraphics[trim=250 50 250 40, clip, width=0.45\textwidth]{figures/fast1550u.jpg}
\includegraphics[width=0.95\textwidth]{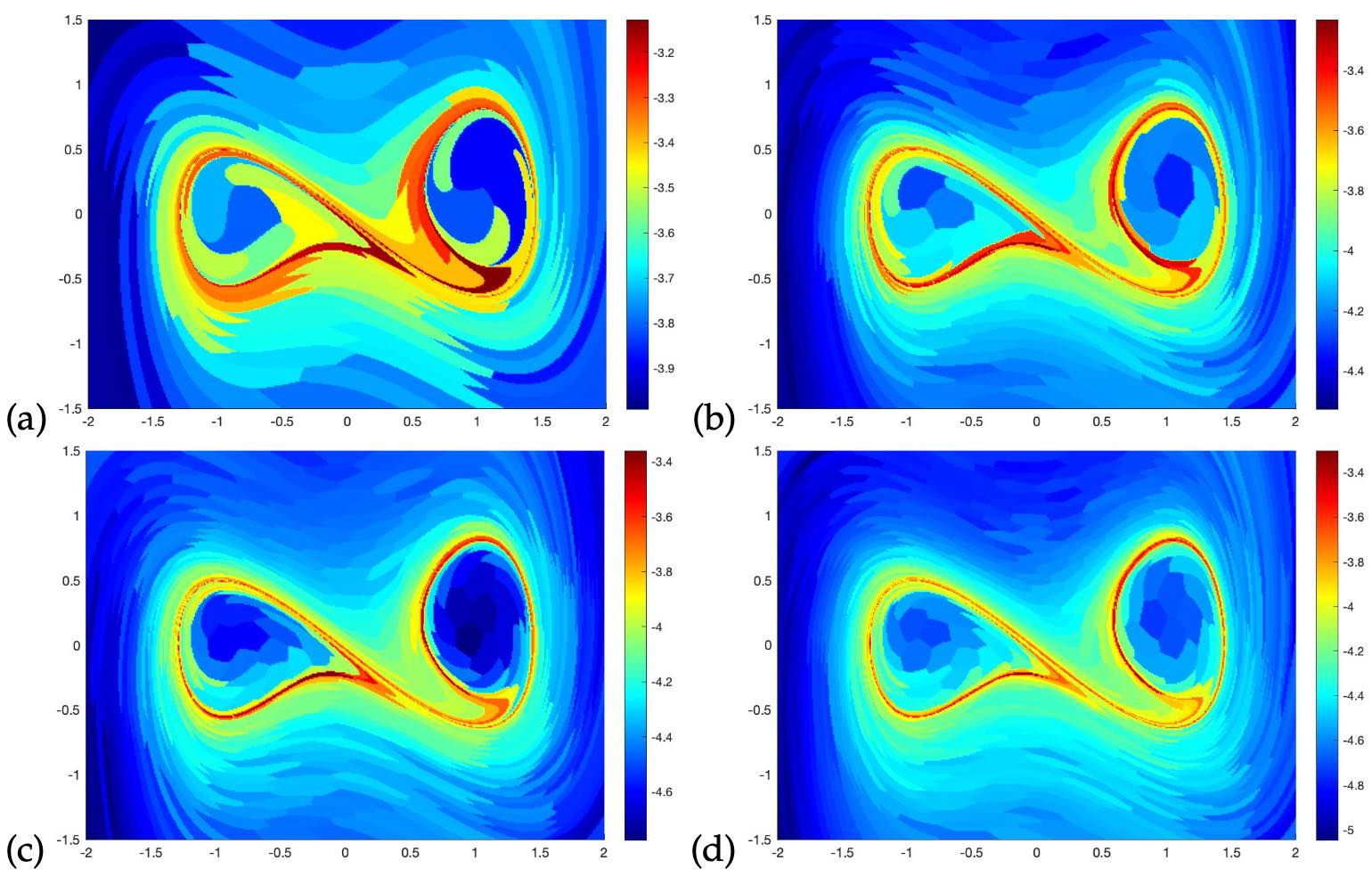}
\caption{(Section \ref{SubSec:Duffing van del Pol}) The WCVE computed by an adaptive refinement approach ($N$=3) with mesh size $\Delta x = \Delta y$ = 1/100 and $\Delta t$ = 0.1 for $T$=10 with (a) $k=$ 100, (b) $k=$ 250, (c) $k=$ 400, and (d) $k=$ 550.}
\label{Ex:Duffing_adaptive}
\end{figure}

\begin{figure}[!htb]
%(a)\includegraphics[trim=250 50 250 40, clip, width=0.45\textwidth]{figures/250umid1.jpg}
%(b)\includegraphics[trim=250 50 250 40, clip, width=0.45\textwidth]{figures/250umid3.jpg}
\includegraphics[width=0.95\textwidth]{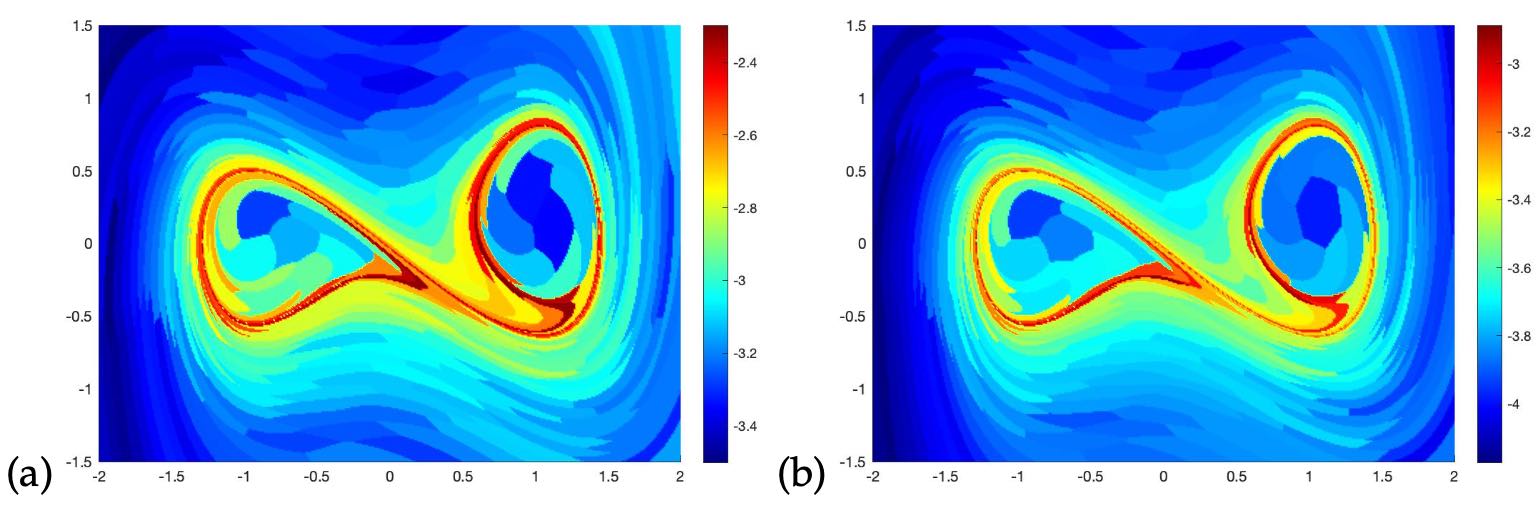}
\caption{(Section \ref{SubSec:Duffing van del Pol}) The intermediate WCVE computed by an adaptive refinement approach with $k=250$ for (a) $N=$3 and (b) $N=$1.}
\label{Ex:Duffing_adaptive_inter_250}
\end{figure}

\begin{figure}[!htb]
%(a)\includegraphics[trim=250 50 250 40, clip, width=0.45\textwidth]{figures/400umid1.jpg}
%(b)\includegraphics[trim=250 50 250 40, clip, width=0.45\textwidth]{figures/400umid3.jpg}
\includegraphics[width=0.95\textwidth]{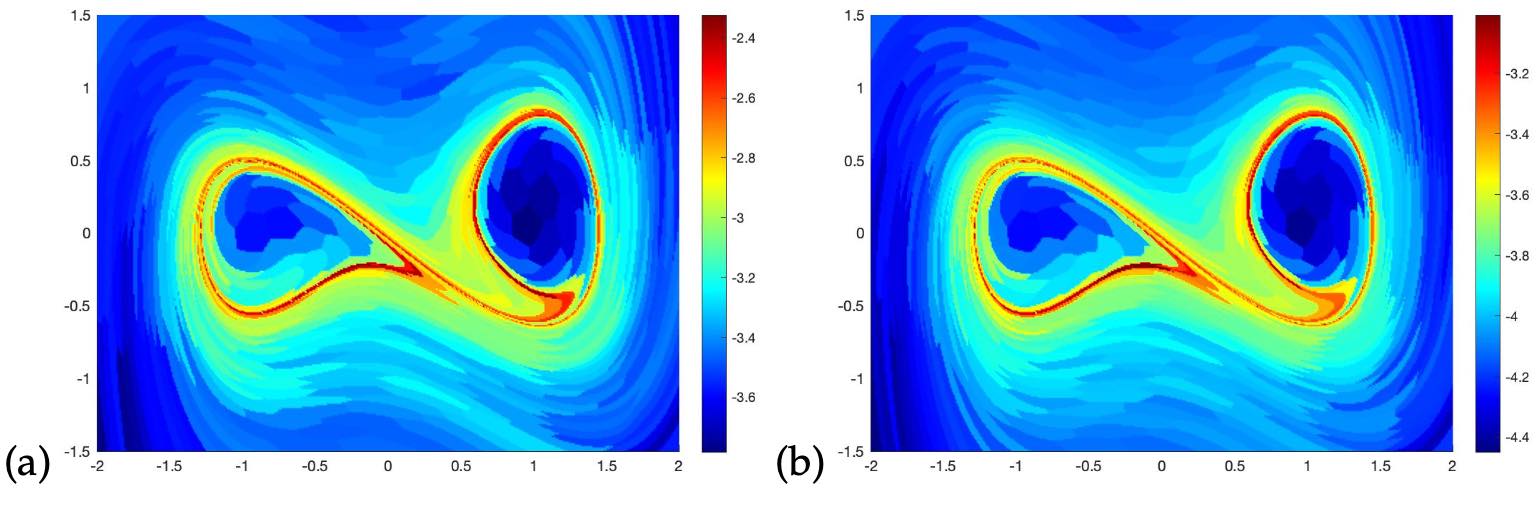}
\caption{(Section \ref{SubSec:Duffing van del Pol}) The intermediate WCVE computed by an adaptive refinement approach with $k=400$ for $N=$3 and (b) $N=$1.}
\label{Ex:Duffing_adaptive_inter_400}
\end{figure}

\begin{table}[!htb]
\centering
(a)
\begin{tabular}{|c|c|c|c|c|c|c|}
\hline
$k$& 100 & 250 & 400 & 550 \\ 
\hline\hline 
Full $k$-means & 488 & 954 & 1355 & 1594 \\ 
Adaptive & 238 & 339 & 406 & 475 \\ 
\reminder{\textit{on-the-fly-5}} & 322 & 618 & 810 & 960\\
\reminder{\textit{on-the-fly-10}} & 289 &535 &702 & 784\\
\hline
\end{tabular} \\
\vspace{0.3cm}
(b)
\begin{tabular}{|c|c|c|c|c|c|c|}
\hline
$k$& 100 & 250 & 400 & 550 \\ 
\hline\hline
Full $k$-means & 807096 & 349202 & 225151 & 167142 \\ 
Adaptive & 802489 & 346400 & 221251 & 162542 \\ 
\reminder{\textit{on-the-fly-5}} & 799819 & 347880 &  223878 & 165452 \\
\reminder{\textit{on-the-fly-10}} & 799087 & 347807 & 224343 & 165410\\
\hline 
\end{tabular}
\caption{(Section \ref{SubSec:Duffing van del Pol}) Comparisons of the full $k$-means algorithm and the adaptive refinement algorithm with measurements taken from the average of 50 trials. (a) The average CPU time (seconds) and (b) the averaged minimized WCSS  (\ref{Eqn:WCSS}). We have rounded all numbers to the closest integer value.}
\label{Ex:Duffing_TimeError}
\end{table}

\begin{figure}[!htb]
%(a)\includegraphics[trim=250 50 250 40, clip, width=0.45\textwidth]{figures/250ut5.jpg}
%(b)\includegraphics[trim=250 50 250 40, clip, width=0.45\textwidth]{figures/250ut52.jpg}
%(c)\includegraphics[trim=250 50 250 40, clip, width=0.45\textwidth]{figures/250ut53.jpg}
%(d)\includegraphics[trim=250 50 250 40, clip, width=0.45\textwidth]{figures/250ut54.jpg}
\includegraphics[width=0.95\textwidth]{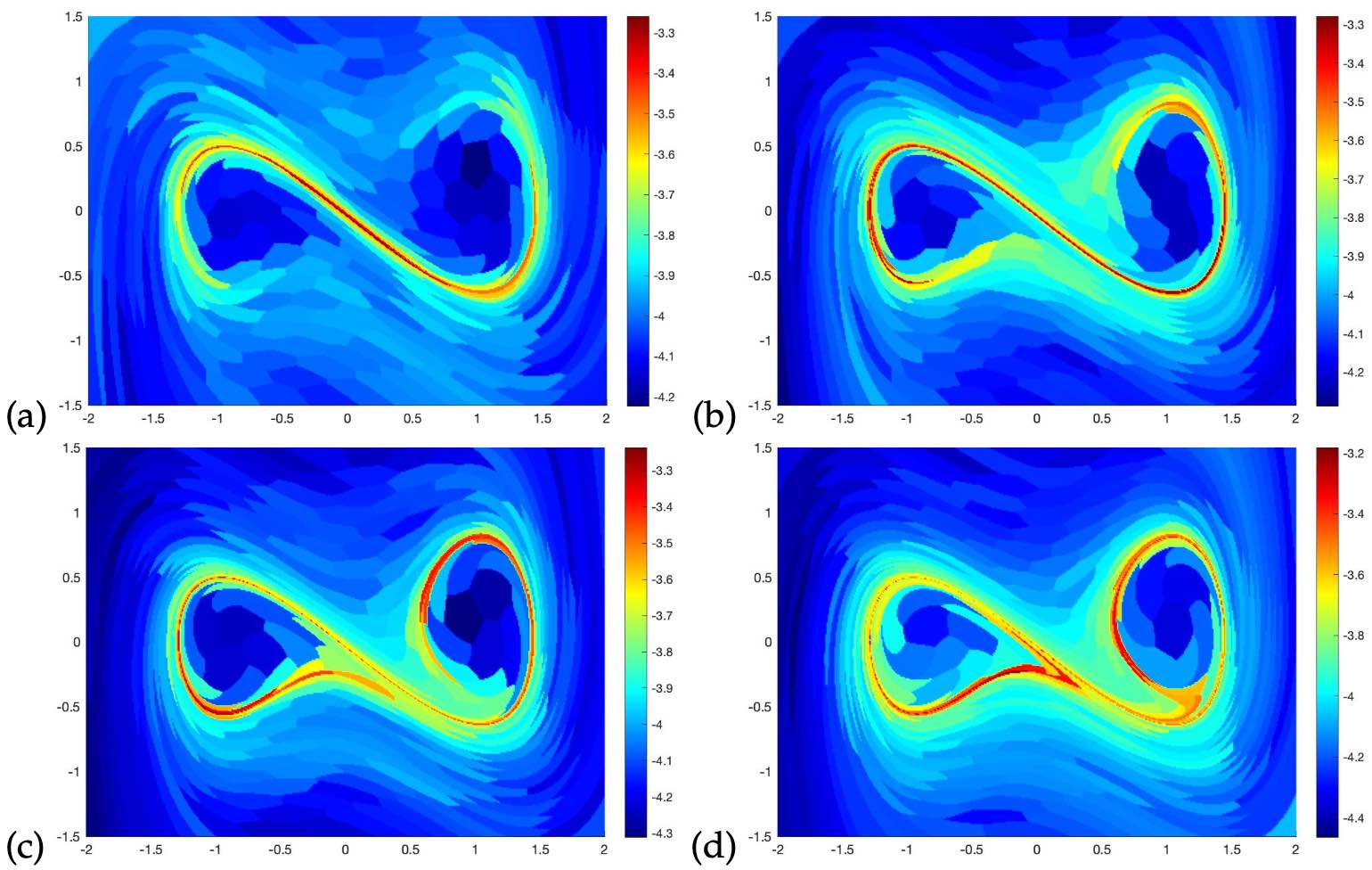}
\caption{(Section \ref{SubSec:Duffing van del Pol}) The WCVE computed by \textit{on-the-fly} approach with $k=250$ for (a) $t_z$= 5, (b) $t_z$= 6.5, (c) $t_z$= 8 and (d) $t_z$= 10. }
\label{Ex:Duffing_fly}
\end{figure}

This example considers the vector field given by 
 \begin{align*}
\mathbf{u}(x,y,t)=\begin{pmatrix} y \\ x-x^3+0.5y(1-x^2)+0.1\sin(t)\end{pmatrix}
\end{align*}
The computational domain for this example is $[-2,2]\times[-1.5,1.5]$, with mesh size $\Delta x = \Delta y$ = 1/100 and time step $\triangle t=0.1$ from $t=0$ to $t=10$. As a reference, we have shown in Figure \ref{Ex:Duffing1} the FTLE computed using the approach as discussed in \cite{youwonleu17}.

Figure \ref{Ex:Duffing_kmean} and Figure \ref{Ex:Duffing_adaptive} demonstrate the obvious ridge of the WCVE field in this complicated dynamics with a relatively coarse mesh size. As the number of $k$ increases, the solutions of the $k$-means algorithm and the adaptive refinement approach can still show clear ridges of the WCVE field. Nevertheless, we see that the overall structure of the WCVE field is in line with the FTLE field in Figure \ref{Ex:Duffing1} for all settings.

The intermediate solutions of an adaptive refinement approach have shown in Figure \ref{Ex:Duffing_adaptive_inter_250} and Figure \ref{Ex:Duffing_adaptive_inter_400}. In general, one can still get a rough structure of the WCVE field with the subsamples of the trajectory for this complicated dynamics system. The ridges of WCVE can be more sharpened with the increase in the time dimensions.

One can improve the computational time and optimization performance significantly for this complicated dynamic system with the adaptive refinement approach, as shown in Table \ref{Ex:Duffing_TimeError}. We see that the CPU time by the adaptive algorithm is approximately one-third of that by the full $k$-mean approach, and the minimized WCSS  from the adaptive algorithm is also smaller than that by the full clustering method. Figure \ref{Ex:Duffing_fly} demonstrates the smooth transformation of the structure of WCVE from $t_z=5$ to $t_z=10$ in the \textit{on-the-fly} approach. 
\reminder{Compared with the full $k$-means clustering algorithm, we find that the CPU time by the \textit{on-the-fly} algorithm is significantly improved as k increases. The minimized WCSS from the \textit{on-the-fly} approach is also smaller than the WCSS from the full clustering method.}

\section*{Acknowledgment}
The second author's work was supported by the Hong Kong RGC under grant 16302819.

\bibliographystyle{plain}
\bibliography{gyou,syleung}

\begin{thebibliography}{10}

\bibitem{bankol17}
R.~Banisch and P.~Koltai.
\newblock Understanding the geometry of transport: diffusion maps for
  lagrangian trajectory data unravel coherent sets.
\newblock {\em Chaos: An Interdisciplinary Journal of Nonlinear Science},
  27(3):035804, 2017.

\bibitem{fkss13}
N.~Ferreira, J.T. Klosowski, C.E. Scheidegger, and C.T. Silva.
\newblock {Vector Field k-Means: Clustering Trajectories by Fitting Multiple
  Vector Fields}.
\newblock {\em Eurographics Conference on Visualization (EuroVis)},
  32(3):201--210, 2013.

\bibitem{fropad15}
G.~Froyland and K.~Padberg-Gehle.
\newblock A rough-and-ready cluster-based approach for extracting finite-time
  coherent sets from sparse and incomplete trajectory data.
\newblock {\em Chaos: An Interdisciplinary Journal of Nonlinear Science},
  25(8):087406, 2015.

\bibitem{grscg07}
S.~Gaffney, A.W. Robertson, P.~Smyth, S.J. Camargo, and M.~Ghil.
\newblock Probabilistic clustering of extratropical cyclones using regression
  mixture models.
\newblock {\em Clim. Dyn.}, 29:423--440, 2007.

\bibitem{gafsmy99}
S.~Gaffney and P.~Smyth.
\newblock {Trajectory Clustering with Mixtures of Regression Models}.
\newblock {\em Proc. 5th ACM SIGKDD Int'l Conf. on Knowledge Discovery and Data
  Mining}, pages 63--72, 1999.

\bibitem{hkth16}
A.~Hadjighasem, D.~Karrasch, H.~Teramoto, and G.~Haller.
\newblock Spectral-clustering approach to lagrangian vortex detection.
\newblock {\em Phys. Rev. E}, 93(6):063107, 2016.

\bibitem{hal01}
G.~Haller.
\newblock Distinguished material surfaces and coherent structures in
  three-dimensional fluid flows.
\newblock {\em Physica D}, 149:248--277, 2001.

\bibitem{hal01b}
G.~Haller.
\newblock {L}agrangian structures and the rate of strain in a partition of
  two-dimensional turbulence.
\newblock {\em Phys. Fluids A}, 13:3368--3385, 2001.

\bibitem{halyua00}
G.~Haller and G.~Yuan.
\newblock {L}agrangian coherent structures and mixing in two-dimensional
  turbulence.
\newblock {\em Physica D}, 147:352--370, 2000.

\bibitem{leehanwha07}
J.~Lee, J.~Han, and K.~Whang.
\newblock Trajectory clustering: a partition-and-group framework.
\newblock {\em ACM SIGMOD International Conference on Management of Data},
  pages 593--604, 2007.

\bibitem{lekshamar07}
F.~Lekien, S.C. Shadden, and J.E. Marsden.
\newblock {L}agrangian coherent structures in $n$-dimensional systems.
\newblock {\em Journal of Mathematical Physics}, 48:065404, 2007.

\bibitem{leu11}
S.~Leung.
\newblock An {E}ulerian approach for computing the finite time {L}yapunov
  exponent.
\newblock {\em J. Comput. Phys.}, 230:3500--3524, 2011.

\bibitem{leu13}
S.~Leung.
\newblock The backward phase flow method for the finite time {L}yapunov
  exponent.
\newblock {\em Chaos}, 23(043132), 2013.

\bibitem{lywn19}
S.~Leung, G.~You, T.~Wong, and Y.K. Ng.
\newblock Recent developments in {E}ulerian approaches for visualizing
  continuous dynamical systems.
\newblock {\em Proceedings of the Seventh International Congress of Chinese
  Mathematicians}, (2):579--622, 2019.

\bibitem{ngleu19}
Y.K. Ng and S.~Leung.
\newblock {Estimating the Finite Time Lyapunov Exponent from Sparse Lagrangian
  Trajectories}.
\newblock {\em Commun. Comput. Phys.}, 26(4):1143--1177, 2019.

\bibitem{oshfed03}
S.~J. Osher and R.~P. Fedkiw.
\newblock {\em Level Set Methods and Dynamic Implicit Surfaces}.
\newblock Springer-Verlag, New York, 2003.

\bibitem{oshset88}
S.~J. Osher and J.~A. Sethian.
\newblock Fronts propagating with curvature dependent speed: algorithms based
  on {Hamilton-Jacobi} formulations.
\newblock {\em J. Comput. Phys.}, 79:12--49, 1988.

\bibitem{shalekmar05}
S.C. Shadden, F.~Lekien, and J.E. Marsden.
\newblock Definition and properties of {L}agrangian coherent structures from
  finite-time {L}yapunov exponents in two-dimensional aperiodic flows.
\newblock {\em Physica D}, 212:271--304, 2005.

\bibitem{vlagunkol02}
M.~Vlachos, D.~Gunopulos, and G.~Kollios.
\newblock {Discovering Similar Multidimensional Trajectories}.
\newblock {\em Proc. 18th Int'l Conf. on Data Engineering}, pages 673--684,
  2002.

\bibitem{wycm11}
J.~Wei, H.~Yu, J.~Chen, and K.~Ma.
\newblock Parallel clustering for visualizing large scientific line data.
\newblock {\em IEEE Symposium on Large Data Analysis and Visualization},
  1:47--55, 2011.

\bibitem{yzzhb17}
D.~Yao, C.~Zhang, Z.~Zhu, J.~Huang, and J.~Bi.
\newblock {Discovering Similar Multidimensional Trajectories}.
\newblock {\em 2017 International Joint Conference on Neural Networks (IJCNN)},
  pages 3880--3887, 2017.

\bibitem{youleu14}
G.~You and S.~Leung.
\newblock An {E}ulerian method for computing the coherent ergodic partition of
  continuous dynamical systems.
\newblock {\em J. Comp. Phys.}, 264:112--132, 2014.

\bibitem{youleu18}
G.~You and S.~Leung.
\newblock Eulerian based interpolation schemes for flow map construction and
  line integral computation with applications to coherent structures
  extraction.
\newblock {\em J. Sci. Comput.}, 74(1):70--96, 2018.

\bibitem{youleu18b}
G.~You and S.~Leung.
\newblock An improved eulerian approach for the finite time lyapunov exponent.
\newblock {\em J. Sci. Comput.}, 76(3):1407--1435, 2018.

\bibitem{youleu20}
G.~You and S.~Leung.
\newblock {Fast Construction of Forward Flow Maps using Eulerian Based
  Interpolation Schemes}.
\newblock {\em J. Sci. Comput.}, 82(32), 2020.

\bibitem{youleu21}
G.~You and S.~Leung.
\newblock Computing the finite time lyapunov exponent for flows with
  uncertainties.
\newblock {\em J. Comput. Phys.}, 425(109905), 2021.

\bibitem{youleu21b}
G.~You and S.~Leung.
\newblock Eulerian algorithms for computing some lagrangian flow network
  quantities.
\newblock {\em J. Comput. Phys.}, 445(11020), 2021.

\bibitem{youwonleu17}
G.~You, T.~Wong, and S.~Leung.
\newblock Eulerian methods for visualizating continuous dynamical systems using
  {L}yapunov exponents.
\newblock {\em SIAM J. Sci. Comp.}, 39(2):A415--A437, 2017.

\end{thebibliography}

\end{document}